\documentclass[12pt]{article}
\usepackage[margin=2.5cm]{geometry}

\usepackage{cite}
\usepackage{amsmath,amssymb,amsfonts}
\usepackage{graphicx}
\usepackage{textcomp}
\usepackage[normalem]{ulem}

\usepackage{longtable}
\usepackage{subfigure}
\usepackage{amsmath}
\allowdisplaybreaks
\usepackage{color}
\usepackage{footnote}
\usepackage{algorithm}
\usepackage{algpseudocode}
\usepackage{algorithmicx}
\usepackage{multirow} 
\usepackage{booktabs}
\usepackage{graphicx}
\usepackage{amssymb}
\usepackage{amsbsy}
\usepackage{array}
\usepackage{longtable}
\usepackage{epstopdf}
\usepackage{pbox}
\usepackage{breqn}
\usepackage{mathrsfs}
\usepackage{multicol}
\usepackage{supertabular}
\usepackage{enumerate}
\usepackage{url}
\usepackage[justification=centering]{caption}
\usepackage{tabu}

\usepackage{enumitem}

\setlist[enumerate]{noitemsep, topsep=0pt}

\newtheorem{thm}{Theorem}

\newtheorem{lmm}{Lemma}

\newtheorem{cor}{Corollary}

\newcommand{\real}{\mathbb{R}} 

\newcommand{\Q}{\mathcal{Q}}

\newcommand{\N}{\mathcal{N}}

\usepackage{etoolbox}
\usepackage{lipsum}

\makeatletter
\let\@xp\expandafter
\DeclareRobustCommand{\qed}{%
  \ifmmode \mathqed
  \else
    \leavevmode\unskip\penalty9999 \hbox{}\nobreak\hfill
    \quad\hbox{\qedsymbol}%
  \fi
}
\let\QED@stack\@empty
\let\qed@elt\relax
\newcommand{\pushQED}[1]{%
  \toks@{\qed@elt{#1}}\@temptokena\expandafter{\QED@stack}%
  \xdef\QED@stack{\the\toks@\the\@temptokena}%
}
\newcommand{\popQED}{%
  \begingroup\let\qed@elt\popQED@elt \QED@stack\relax\relax\endgroup
}
\def\popQED@elt#1#2\relax{#1\gdef\QED@stack{#2}}
\newcommand{\qedhere}{%
  \begingroup \let\mathqed\math@qedhere
    \let\qed@elt\setQED@elt \QED@stack\relax\relax \endgroup
}
\newif\ifmeasuring@
\newif\iffirstchoice@ \firstchoice@true
\def\setQED@elt#1#2\relax{%
  \ifmeasuring@
  \else \iffirstchoice@ \gdef\QED@stack{\qed@elt{}#2}\fi
  \fi
  #1%
}
\def\qed@warning{%
  \PackageWarning{amsthm}{The \@nx\qedhere command may not work
    correctly here}%
}
\newcommand{\mathqed}{\quad\hbox{\qedsymbol}}
\def\linebox@qed{\hfil\hbox{\qedsymbol}\hfilneg}
\@ifpackageloaded{amsmath}{%
  \def\math@qedhere{%
    \@ifundefined{\@currenvir @qed}{%
      \qed@warning\quad\hbox{\qedsymbol}%
    }{%
      \@xp\aftergroup\csname\@currenvir @qed\endcsname
    }%
  }
  \def\displaymath@qed{%
    \relax
    \ifmmode
      \ifinner \aftergroup\linebox@qed
      \else
        \eqno
        \let\eqno\relax \let\leqno\relax \let\veqno\relax
        \hbox{\qedsymbol}%
      \fi
    \else
       \aftergroup\linebox@qed
    \fi
  }
  \@xp\let\csname equation*@qed\endcsname\displaymath@qed
  \def\equation@qed{%
    \iftagsleft@
      \hbox{\phantom{\quad\qedsymbol}}%
      \gdef\alt@tag{%
        \rlap{\hbox to\displaywidth{\hfil\qedsymbol}}%
        \global\let\alt@tag\@empty
      }%
    \else
      \gdef\alt@tag{%
        \global\let\alt@tag\@empty
        \vtop{\ialign{\hfil####\cr
                \tagform@\theequation\cr
                \qedsymbol\cr}}%
        \setbox\z@
      }%
    \fi
  }
  \def\qed@tag{%
    \global\tag@true \nonumber
    &\omit\setboxz@h {\strut@ \qedsymbol}\tagsleft@false
    \place@tag@gather
    \kern-\tabskip
    \ifst@rred \else \global\@eqnswtrue \fi \global\advance\row@\@ne \cr
  }
  \def\split@qed{%
    \def\endsplit{\crcr\egroup \egroup \ctagsplit@false \rendsplit@
      \aftergroup\align@qed
    }%
  }
  \def\align@qed{%
    \ifmeasuring@ \tag*{\qedsymbol}%
    \else \let\math@cr@@@\qed@tag
    \fi
  }
  \@xp\let\csname align*@qed\endcsname\align@qed
  \@xp\let\csname gather*@qed\endcsname\align@qed
}{
  \def\math@qedhere{%
    \@ifundefined{\@currenvir @qed}{%
      \qed@warning \aftergroup\displaymath@qed
    }{%
      \@xp\aftergroup\csname\@currenvir @qed\endcsname
    }%
  }
  \def\displaymath@qed{%
    \relax
    \ifmmode
      \ifinner \aftergroup\aftergroup\aftergroup\linebox@qed
      \else
        \eqno \def\@badmath{$$}%
        \let\eqno\relax \let\leqno\relax \let\veqno\relax
        \hbox{\qedsymbol}%
      \fi
    \else
       \aftergroup\linebox@qed
    \fi
  }
  \@ifundefined{ver@leqno.clo}{%
    \def\equation@qed{\displaymath@qed \quad}%
  }{%
    \def\equation@qed{\displaymath@qed}%
  }
  \def\@tempa#1$#2#3\@nil{%
    \def\[{#1$#2\def\@currenvir{displaymath}#3}%
  }%
  \expandafter\@tempa\[\@nil
}
\@ifpackageloaded{amstex}{%
  \def\@tempa{TT}%
}{%
  \@ifpackageloaded{amsmath}{%
    \def\@tempb#1 v#2.#3\@nil{#2}%
    \ifnum\@xp\@xp\@xp\@tempb\csname ver@amsmath.sty\endcsname v0.0\@nil
       <\tw@
      \def\@tempa{TT}%
    \else
      \def\@tempa{TF}%
    \fi
  }{%
    \def\@tempa{TF}
  }%
}
\if\@tempa
  \renewcommand{\math@qedhere}{\quad\hbox{\qedsymbol}}%
\fi
\newcommand{\openbox}{\leavevmode
  \hbox to.77778em{%
  \hfil\vrule
  \vbox to.675em{\hrule width.6em\vfil\hrule}%
  \vrule\hfil}}
\DeclareRobustCommand{\textsquare}{%
  \begingroup \usefont{U}{msa}{m}{n}\thr@@\endgroup
}
\providecommand{\qedsymbol}{\openbox}
\makeatother

\newtheorem{theorem}{Theorem}
\preto{\theorem}{\pushQED{\qed}}
\preto{\endtheorem}{\popQED}

\begin{document}
\title{Resilient Ramp Control for Highways Facing Stochastic Perturbations}
\author{
Yu Tang,
Li Jin,
Alexander A. Kurzhanskiy
and Saurabh Amin
\thanks{This work was in part supported by NYU Tandon School of Engineering, C2SMART University Transportation Center, US NSF Award CMMI-1949710 and CAREER Award CNS-1453126, SJTU UM Joint Institute, J. Wu \& J. Sun Endowment Fund, and Singapore NRF Future Urban Mobility.}
\thanks{L. Jin is with the UM Joint Institute and the School of Electronic Information and Electrical Engineering, Shanghai Jiao Tong University, China. L. Jin is also and Y. Tang is with the Tandon School of Engineering, New York University, USA.
A. A, Kurzhanskiy is with the Institute of Transportation Studies, University of California, Berkeley, USA.
S. Amin is with the Laboratory of Information and Decision Systems and the Department of Civil and Environmental Engineering, Massachusetts Institute of Technology, USA
(emails: li.jin@sjtu.edu.cn,  tangyu@nyu.edu, akurzhan@berkeley.edu, amins@mit.edu).}
}
\maketitle

\begin{abstract}                          

Highway capacity is often subject to stochastic perturbations due to the combined effects of weather, traffic mixture, driver behavior, etc. This paper is motivated by the need of a systematic approach to traffic control with performance guarantees in the face of such perturbations. We develop a novel control-theoretic method for designing perturbation-resilient ramp metering. We consider a cell-transmission model with 1) Markovian cell capacities and 2) buffers representing on-ramps and upstream mainline. Using this model, we analyze the stability of 
on-ramp queues by constructing piecewise Lyapunov functions that consider the nature of nonlinear traffic dynamics. Then, we design ramp controllers that guarantee  bounds for throughput and queue sizes. We also formulate the problem of coordinated ramp metering as a bi-level optimization with non-convex inner sub-problems. To address the computational issue in solving this problem, we also consider localized and partially coordinated reformulations. A case study of a 18.1-km highway in Los Angeles, USA indicates a 8.3\% (resp. 9.9\%) reduction of vehicle-hours-traveled obtained by the localized (resp. partially coordinated) control, both outperforming the classical ALINEA and METALINE controllers.
\end{abstract}


\section{Introduction}

\subsection{Motivation}
Modern control systems continue to play an important role in improving highway system performance. Technologies such as embedded sensors and dynamic traffic signal controllers have enabled several highway control strategies, including ramp metering, variable speed limit, and dynamic routing \cite{papageorgiou2003review,kurzhanskiy10,ferrara2018freeway}. Most of these strategies are designed for nominal (and deterministic) settings. However, the intrinsic uncertainties faced by highway systems can compromise the performance of such strategies  \cite{jin14, jin19tac}. In this paper, we focus on control of highway on-ramps under capacity variations \cite{polus02,banks1991two}\footnote{We use the terms capacity variations/perturbations/fluctuations interchangeably.}. 

Usually, \emph{highway capacity} is defined as the maximum sustained flow rate under prevailing roadway and traffic conditions. In practice, the capacity may change when the traffic is perturbed by weather \cite{heshami2019deterministic}, driving behavior \cite{khoshyaran2015capacity}, traffic merging \cite{asgharzadeh2020effect}, and so on. Besides, such perturbations are often hard to predict and prevent. To address such issues, one can consider a stochastic model of highway capacity \cite{polus02}. The following example shows that highway dynamics can be better described by modeling stochastic capacities.
\begin{figure}[htbp]
\centering
\subfigure[A section of Interstate 210 Eastbound in Los Angeles, California, U.S..]{
\centering
\includegraphics[width=0.65\linewidth]{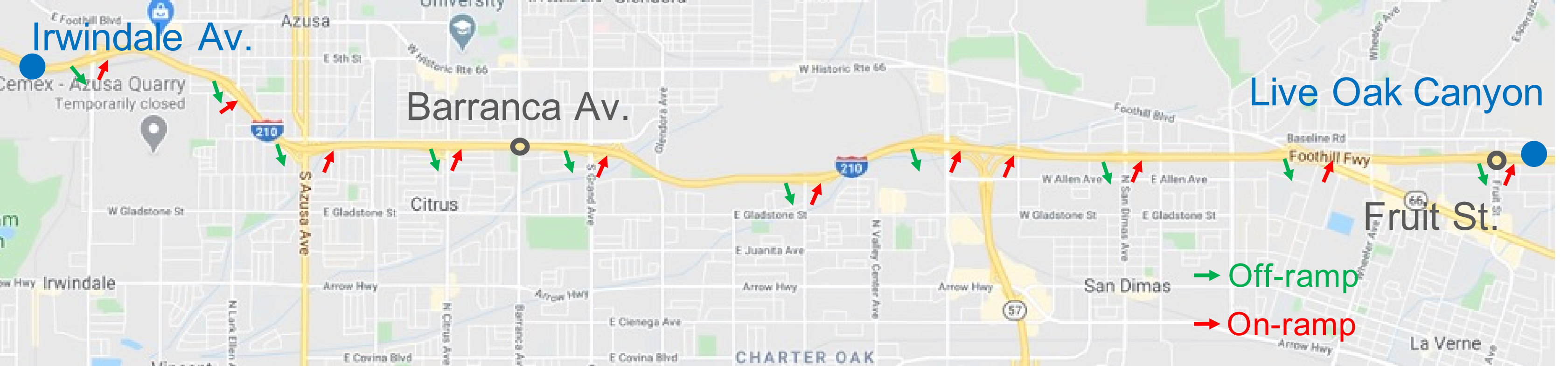}
\label{fig_map}
}

\subfigure[Barranca Avenue.]{
\centering
\includegraphics[width=0.25\linewidth]{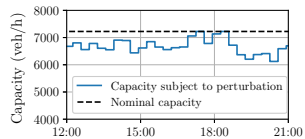}
\label{fig_cap_cell7}
}
\subfigure[Fruit Street.]{
\centering
\includegraphics[width=0.25\linewidth]{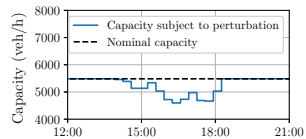}
\label{fig_cap_cell16}
}
\caption{Capacity fluctuations at two bottlenecks.}
\label{fig_cap}
\end{figure}
\begin{figure*}[h!]
\centering
\subfigure[True data.]{
\centering
\includegraphics[width=0.25\textwidth]{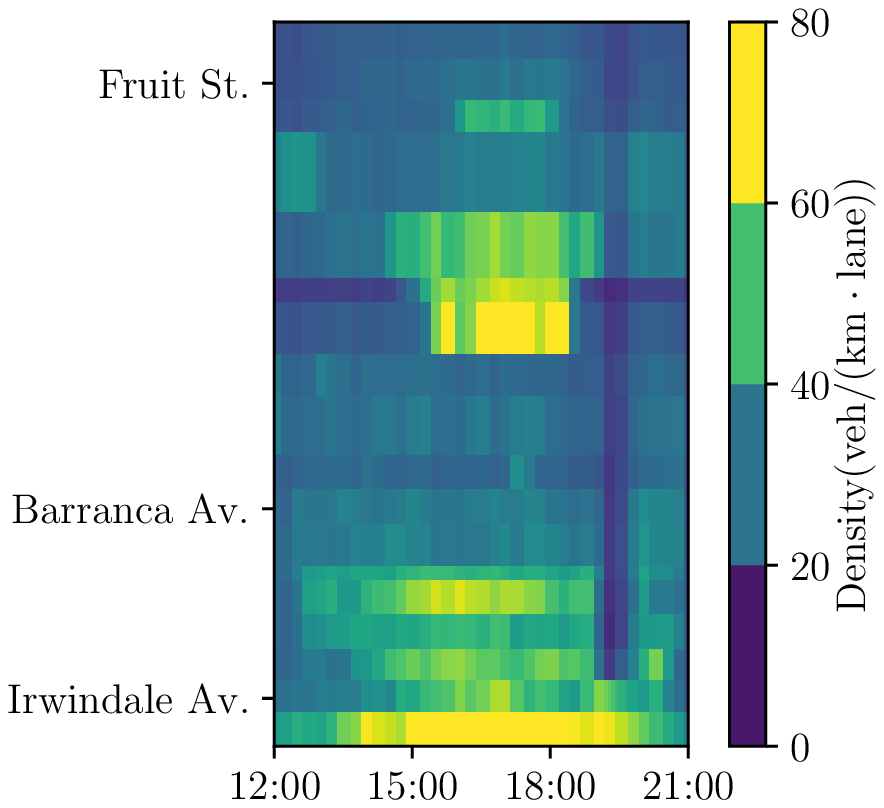}
\label{fig_truespeedmap}
}
\subfigure[CTM.]{
\centering
\includegraphics[width=0.25\textwidth]{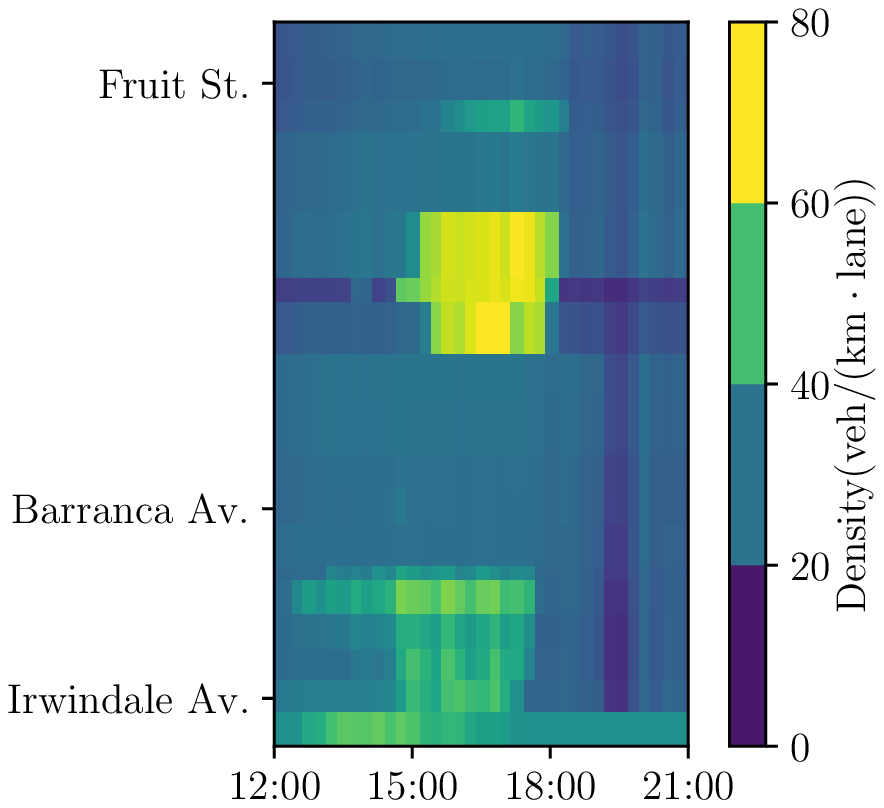}
\label{fig_ctmspeedmap}
}
\subfigure[SS-CTM.]{
\centering
\includegraphics[width=0.25\textwidth]{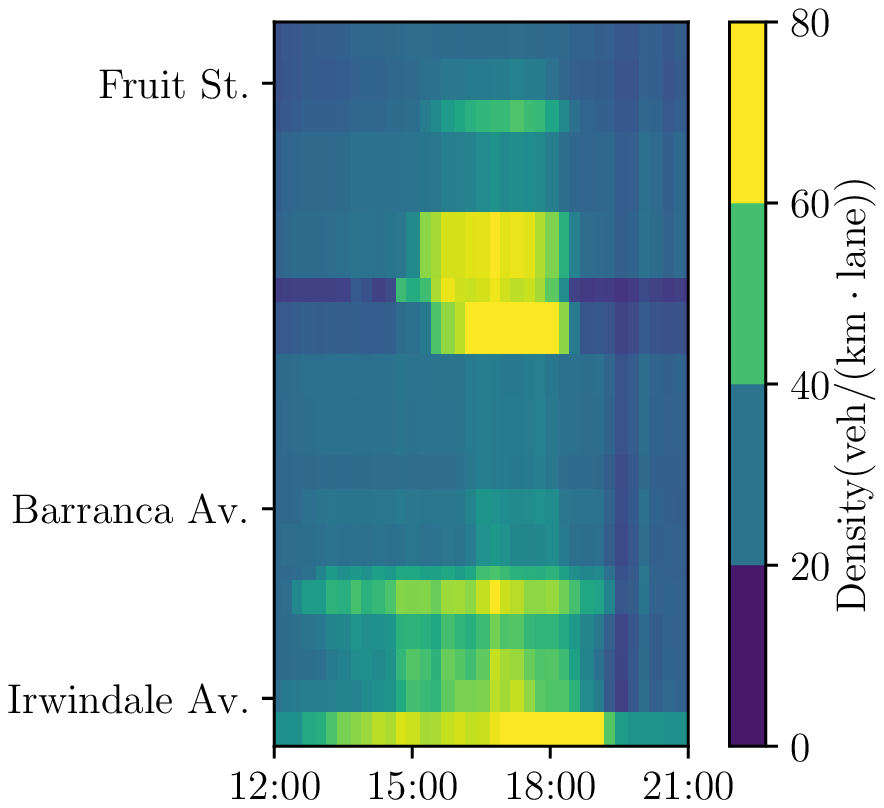}
\label{fig_ssctmspeedmap}
}
\caption{Comparison of traffic density maps.}
\label{fig_calibratedSpeedMap}
\end{figure*}

We consider a 18.1-km Interstate 210 Eastbound (I-210E) section as shown in Fig.~\ref{fig_map}. We calibrated both a conventional cell-transmission model (CTM \cite{daganzo94}) and a stochastic-switching cell transmission model (SS-CTM \cite{jin19tac}) with Markovian capacities. In SS-CTM, the capacities switch between different states according to a Markov process. This approach has been reported as a promising one for modeling stochastic traffic flows \cite{sumalee2011stochastic}, occurrence of traffic breakdowns \cite{evans2001probability}, and vehicle headway/spacing \cite{chen2010markov}. We calibrated this model based on \cite{jin2017calibration}; see related methods in \cite{polus02}. Our calibration uses PeMS data \cite{varaiya09} (see Appendix~A). Figs.~\ref{fig_cap_cell7} and \ref{fig_cap_cell16} show capacity fluctuation at traffic bottlenecks near Barranca Ave. and Fruit St. in Los Angeles on March 26, 2019.

We also simulated and compared these models in terms of mean absolute percentage error (MAPE) of vehicle travel times.  The deterministic CTM gives a MAPE 8.9\%, while the SS-CTM yields a MAPE of 4.3\%. Fig.~\ref{fig_calibratedSpeedMap} visualizes the real and simulated traffic density maps. The CTM underestimates the traffic jam between Irwindale Ave. and Barranca Ave., while the SS-CTM significantly reduces this modeling error by introducing stochastic capacities. This example motivates our hypothesis that control design based on SS-CTM may be more promising than CTM-based control in managing congestion resulting from capacity variations.

In this paper, we develop a control-theoretic approach to designing ramp control that is resilient to Markovian capacity. The proposed approach can guarantee a lower bound for throughput and an upper bound for on-ramp queue sizes in the face of perturbations. To this end, we consider SS-CTM, a spatially discretized traffic flow model with stochastic capacities. This model extends the well-known CTM by introducing two features: 1) Highway capacities varying according to a Markov process; and 2) Infinite-sized buffers representing on-ramps and upstream mainline. Besides tracking mainline traffic densities analogous to CTM, this model uses a Markov process for time-varying cell capacities and infinite-sized buffers for queues at on-ramps and upstream mainline. Using this model, we study the stability of queues under a class of ramp metering policies. The stability analysis enables us to formulate a bi-level program for designing a coordinated ramp metering strategy (similar to \cite{papageorgiou1990modelling,bhouri2013isolated}). To address the computational difficulty arising from the non-convex and large-scale nature of this formulation, we consider localized and partially coordinated reformulations. The former relies only on local information (similar to \cite{papageorgiou91,wang2014local}) and does not consider the interaction between multiple ramps;  the latter sequentially designs ramp controllers such that the upstream controllers depend on the downstream ones. We apply both approaches to the above I-210E highway system and demonstrate that they lead to higher performance, in comparison to those assume deterministic highway capacity.

\subsection{Related work}
Prior literature has investigated various ramp control schemes. Papageorgiou et al. proposed the integral controller ALINEA for localized ramp metering \cite{papageorgiou91} and then generalized it to the coordinated METALINE \cite{papageorgiou1990modelling}. These algorithms were validated in a series of field tests \cite{papageorgiou1997alinea} and prompted a range of controllers \cite{wang2014local,frejo2018feed,pasquale2020hierarchical}. These have significantly refined the original ALINEA; however they still rely on nominal critical density that is defined based on maximal capacity. Another way is to formulate the on-ramp control as an optimization problem based on deterministic traffic flow models. Optimal ramp metering over a finite horizon was addressed via nonlinear optimization \cite{kotsialos2004nonlinear,gomes06,schmitt2018exact,como2016convexity}. Real-time ramp metering can be also designed using model predictive control (MPC) \cite{hegyi2005model,muralidharan2012optimal}, based on computationally efficient formulations. A distributed MPC-based ramp metering has been also proposed to improve tractability \cite{ferrara2014distributed}. 

An alternative is to design optimal or near-optimal static state feedback controllers with simple structures \cite{jafari2019structural}, which are considered in our paper. Our controllers do not rely on real-time computation and can be designed based on offline-calibrated SS-CTM. Furthermore, they can be flexibly implemented in either a centralized or a decentralized manner. Here we build on our previous work which showed that capacity perturbations may destabilize on-ramp queues that are stable in the deterministic, nominal/average setting \cite{jin19tac}. 
In particular, capacity fluctuations can induce new bottlenecks and the nominally designed controllers can fail to mitigate traffic queues at these bottlenecks. Some attempts to solve this problem have been reported in the literature. One approach is to estimate capacity drop and design control laws based on a modified fundamental diagram \cite{zhang2018stability}. Although the capacity drop in this method is explicitly modeled by the relation between capacity and density, it does not account for perturbations caused by other factors. Robust control may be also applied based on appropriately defined uncertainty sets \cite{como13i,chow2014robust, zhong14,schmitt2020convexity}. However, our paper is expressly driven by the need to handle stochastic capacity perturbations, including the ones with low probability of occurrence. There exists recent work that studies highway control under stochastic perturbations \cite{mehr2016probabilistic, heshami2021ramp}. However, the body of work is somewhat limited in its treatment of stability and does not provide resiliency guarantees against stochastic capacity. 

\subsection{Our contributions}
We formalize the notion of resilient ramp metering for stochastic capacity, and address two main control design questions:
\begin{enumerate}
    \item[(i)] How to formulate an optimization model with stability constraint and control cost related to queue size or throughput?
    \item[(ii)] How to obtain practically relevant control policies by solving this model when ramp controllers are coupled?
\end{enumerate}

We address (i) by exploiting the Foster-Lyapunov criterion \cite{meyn93}. This criterion leads to the \emph{drift condition} for the feedback-controlled SS-CTM and it is sufficient to ensure stability of the controllers. The generic form of Foster-Lyapunov criterion requires verifying the drift condition everywhere over the state space, but we simplify it by constructing a bounded invariant set (a subset of the state space) and show that one can instead verify the drift condition over the invariant set. Second, we utilize the drift condition to obtain an upper bound of the buffer queue size for a given control law and traffic demand, and a corresponding lower bound of throughput. Our design objective can thus be formulated as problem of minimizing the queue size given traffic demand subject to the stability constraint. Importantly, we use the upper bound as a proxy for the exact queuing cost. The numerical examples in Sections~\ref{sec_distributed} and \ref{sec_cordinated} show that this bound is reasonably tight in practice.

To address (ii), we develop localized and partially coordinated reformulations of the general coordinated control design problem. These reformulations may reduce control efficiency somewhat, but provide tractable means to solve large-scale problems. The localized ramp metering is designed independently of other ramp controllers. For this reformulation, we propose a Lyapunov function based on the insights about the traffic dynamics; this simplifies both the stability condition (Theorem~\ref{thm_distributed}) and control objective. Then the controller can be obtained by solving a bi-level program with a non-convex inner problem. 
The technique of localized reformulation can also be applied to the fully coordinated control design (Theorem~\ref{thm_centralized}), which inspires a partially coordinated reformulation. We note that a fully coordinated design accounts for the interactions among ramp controllers, while the partially coordinated approach only considers the impacts of downstream ramp metering on the upstream ones. Essentially, the latter decomposes the coordination problem into a set of interconnected sub-problems and solves them one-by-one from the downstream to the upstream. We also present the corresponding stability condition (Theorem~\ref{thm_cascaded}) and numerically show in Section~\ref{sec_cordinated} that it can yield a good approximate solution.

Finally, we validate our approach using the I-210E case study introduced in Fig.~\ref{fig_map}. Two benchmarks are considered, one with the localized ALINEA and the other with the coordinated METALINE. Numerical evaluation shows that, compared to ALINEA and METALINE, our approach is more efficient in terms of reducing total vehicle hours (both on the mainline and at the on-ramps).
Our localized and partially coordinated control design reduced total vehicle hours by 8.3\% and 9.9\% respectively, while ALINEA and MEATLINE contributed reductions of 5.1\% and 6.2\% respectively. Our results suggest that the existing methods for classical ramp metering can be refined to account for stochastic capacity perturbations, and our approach is useful in providing resiliency guarantees in terms to stability and performance.

The rest of this paper is organized as follows. We first introduce the SS-CTM and formulate the control problem in Section~\ref{sec_model}. The following Sections~\ref{sec_distributed} and \ref{sec_cordinated} focus on localized and coordinated ramp metering, respectively. Then, Section~\ref{sec_simulate} presents our case study. Finally, we summarize the paper and discuss the future research in Section~\ref{sec_conclude}.
\section{Modeling and formulation}
\label{sec_model}


\subsection{Stochastic traffic model}
\begin{figure}[htb]
\centering
\includegraphics[width=0.5\linewidth]{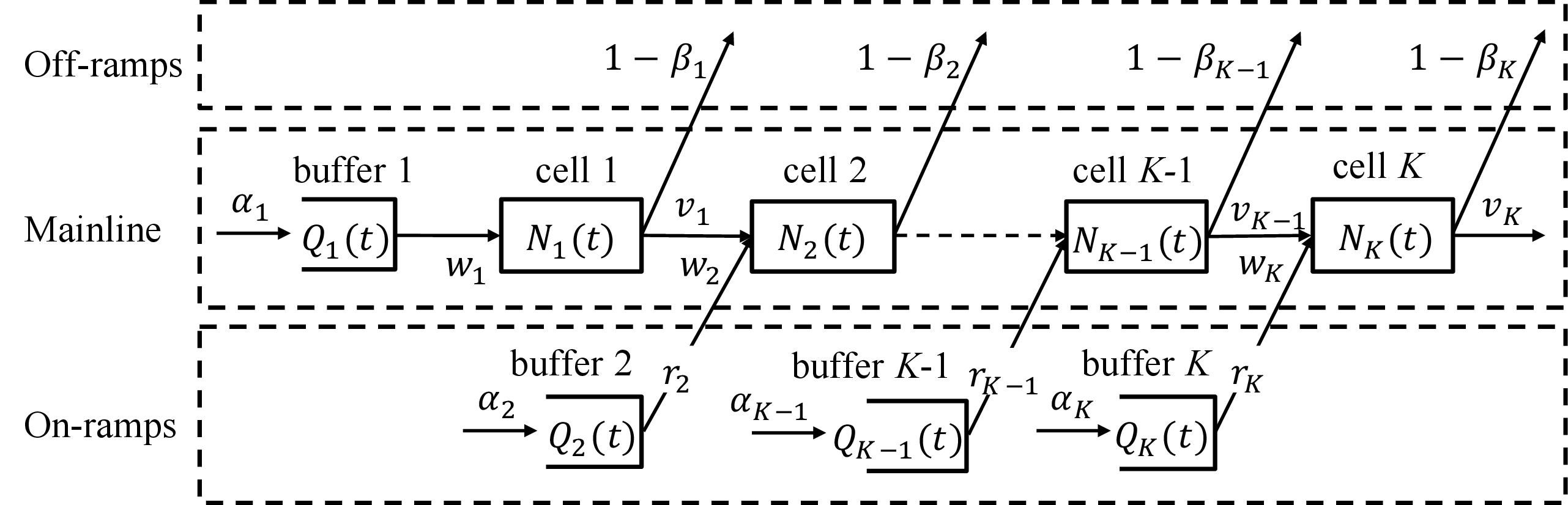}
\caption{A highway with $K$ cells, $K$ buffers and $K$ off-ramps.}
\label{fig_ctm}
\end{figure}

We consider a highway with $K$ mainline \emph{cells}, $K$ on-ramp \emph{buffers}, and $K$ off-ramps, as shown in Fig.~\ref{fig_ctm}. The first buffer is not an actual on-ramp; instead, it represents the upstream highway section and stores the upstream mainline traffic. Each cell $k$ has the following parameters: \emph{length} $l_k$ (km), \emph{free-flow speed} $v_k$ (km/hr), \emph{congestion wave speed} $w_k$ (km/hr), \emph{jam density} $n^{\mathrm{jam}}_k$ (veh/km) and \emph{mainline ratio} $\beta_k\in[0, 1]$. Here $\beta_k$ denotes the fraction of traffic from cell $k$ entering cell $k+1$; the remaining traffic flow leaves the highway at the $k$th off-ramp. Without loss of generality, we let $\beta_K=0$.
The $k$th buffer has a capacity ${U_k}$ (veh/hr) and is subject to a demand $\alpha_k\in[0, U_k]$. We assume that both $\alpha_k$ and $\beta_k$ are constant. The main reason is that the traffic model is considered for design of ramp metering which in practice is only turned on in peak hours \cite{adot} when traffic demands are steady. Note that ramp metering regulates the flow $r_k$ from buffer $k$ to cell $k$ for $k=2,3,\cdots,K$; see Fig.~\ref{fig_ctm}.
The \emph{continuous state} of the highway is $(q,n)$, where $q=[q_1,\cdots,q_K]^{\mathrm{T}}$ is the vector of buffer queues and $n=[n_1,\cdots,n_K]^{\mathrm{T}}$ is the vector of traffic densities in mainline cells.
The set of permissible queue lengths is $\Q:=\real_{\ge0}^K$, and the set of traffic densities is $\N:=\prod_{k=1}^K[0,n^{\mathrm{jam}}_k]$.
Hence, the \emph{continuous state space} is $\mathcal Q\times\mathcal N$.

In our stochastic traffic flow model, the state of the Markov process is referred to as \emph{mode}, and we denote by  $\mathcal S:=\{1,2,\ldots,m\}$ the set of modes. Every mode $s\in\mathcal{S}$ is associated with a vector of cell capacities $[F_{1,s},\cdots, F_{K,s}]^{\mathrm{T}}$. Thus stochastic capacity is captured by the mode transition. Let $S(t)$ be the mode at time $t$. We assume that the Markov process $\{S(t);t\ge0\}$ is ergodic. It indicates that the Markov process returns to any state in a finite time almost surely. This assumption is reasonable, considering recurrent traffic jam day by day and underlying capacity variation \cite{jin19tac}.

Let $\{\lambda_{s,s'};s,s'\in\mathcal S\}$ denote the inter-mode transition rates.
Without loss of generality, we let $\lambda_{s,s}=0$ for all $s\in\mathcal{S}$.
Hence, we have $\Pr\{S(t+\delta)=s'|S(t)=s\}=\lambda_{s,s'}\delta+o(\delta)$, where $\delta$ is an infinitesimal time increment. The assumption of ergodicity means that a unique vector of steady-state probabilities $ p=[p_1, \cdots,p_m]^{\mathrm{T}}$ exists  such that
\begin{subequations}
\begin{align}
    &\sum_{s'\in\mathcal S}p_{s'}\lambda_{s',s}=\sum_{s'\in\mathcal S}p_{s}\lambda_{s,s'}, ~s\in\mathcal S, \label{eq_prob_a}\\
    &\sum_{s\in\mathcal S}p_{s}=1,  \label{eq_prob_b} \\
    &p_s\ge0, ~s\in\mathcal S, \label{eq_prob_c}
\end{align}
\end{subequations}
where \eqref{eq_prob_a} denotes the equilibrium of steady-state distribution. We assume that each mode $s$ associated with $[F_{1,s},\cdots,F_{K,s}]^{\mathrm{T}}$ and the transition rates $\lambda_{s,s'}$ are calibrated in advance. The calibration details are available in \cite{jin2017calibration}.

For each $k$, we define
\begin{align*}
    \bar F_k:=\sum_{s\in\mathcal S}p_sF_{k,s}, ~
    F^\max_k:=\max_{s\in\mathcal S}F_{k,s}, ~F^\min_k:=\min_{s\in\mathcal S}F_{k,s},
\end{align*}
which can be interpreted as the mean, the maximal, and the minimal capacity of cell $k$, respectively. We recognize $F^\max_k$ as the nominal capacity of cell $k$ and define the \emph{nominal critical density} as $n_k^c := F^{\max}_k/v_k$. Following \cite{daganzo94}, we assume
\begin{equation}
v_kn_k^c \leq w_k(n_k^{\mathrm{jam}}-n_k^c), \label{eq_fun}
\end{equation}
which implies the receiving flow is not less than the sending flow at the nominal critical density $n_k^c$.

For notational convenience, we use $\phi=(s,q,n)$ 
to denote hybrid state and $\Phi(t)=(S(t),Q(t),N(t))$
to denote (hybrid) stochastic process. Using this convention, we let $Q_k(t)$ be the queue length in the $k$th on-ramp and ${N}_k(t)$ denote the {traffic density} in the $k$th cell at time $t$, and thus $Q(t)=[Q_1(t),\cdots,Q_K(t)]^{\mathrm{T}}\in\mathcal Q$ and $N(t)=[N_1(t),\cdots,N_K(t)]^{\mathrm{T}}\in\mathcal N$.

The following specifies the flow functions and traffic dynamics. First, the \emph{mainline inflow} $r_1$ is given by
\begin{equation}
    r_1(q_1, n_1) = 
    \begin{cases}
    \min\{\alpha_1, w_1(n_1^{\mathrm{jam}} - n_1)\} & q_1=0, \\
    \min\{U_1, w_1(n_1^{\mathrm{jam}} - n_1)\} & q_1 > 0.
    \end{cases}
\end{equation}
For $k=2,\cdots,K$, the \emph{controlled on-ramp flow} $r_k^\mu$ from buffer $k$ to cell $k$ is given by
\begin{equation}
    r_k^\mu(q_k, n) =
    \begin{cases}
    \min\{\alpha_k, w_k(n_k^{\mathrm{jam}} - n_k),\mu_k(n)\} & q_k=0, \\
    \min\{U_k, w_k(n_k^{\mathrm{jam}} - n_k),\mu_k(n)\} & q_k > 0,
    \end{cases}
\end{equation}
where $\mu_k(n)$ denotes the control input from on-ramp $k$. Then the \emph{outflow} ${f_k^\mu}$ from cell $k$ is given by
\begin{subequations}
\begin{align}
f_k^\mu(\phi) =& \min\Big\{v_kn_k,F_{k,s}, \frac{1}{\beta_k}\Big(w_{k+1}(n^{\mathrm{jam}}_{k+1}-n_{k+1}) \nonumber \\
&\quad\quad - r_{k+1}^\mu(q_{k+1}, n)\Big)\Big\},~k=1,\cdots,K-1, \label{eq_fk} \\
f_K^\mu(\phi) =& \min\{v_Kn_K,F_{K,s}\}. \label{eq_fK}
\end{align}
\end{subequations}
Note that \eqref{eq_fk} indicates higher merging priority of on-ramp flows and the first-in-first-out rule for off-ramp flows \cite{ferrara2018freeway}. Also \eqref{eq_fK} indicates no bottlenecks downstream of cell $K$.


We consider that on-ramp flows are regulated by \emph{affine controllers} \cite{skaf2010design}. This control law is favored because the simple state feedback can achieve performance comparable to that of more sophisticated methods \cite{schmitt2017sufficient}. Specifically, the control policy $\mu:\mathcal N\to\mathbb R_{\ge0}^{K-1}$ is given by
\begin{align}
    \mu_k(n)= u_k- \kappa_k n_k,
    ~k=2,\ldots,K, \label{eq_controller}
\end{align}
where the controller $\mu_k$ parameterized by $u_k$ and $\kappa_k$ limits the ramp flow based on the traffic density $n_k$. It indicates that $\mu_k(n)$ monotonically decreases as the downstream traffic density $n_k$ increases. The proposed ramp metering \eqref{eq_controller} relies on local measurements, similar to the classical ALINEA. 

In the following, we denote by $r_k^\mu(q_k,n_k)$ the on-ramp flows. The control parameters $u_k$ and $\kappa_k$ can be optimized in either a localized or a coordinated manner. The optimal values of these parameters rely on the problem settings, such as traffic demand and capacity fluctuation. For example, when mainline demand increases, the designed ramp metering tends to be more aggressive with smaller $u_k$ or larger $\kappa_k$; see Fig.~\ref{fig_twocell_analysis}. Note that our approach does not necessarily require an affine control. More sophisticated ramp metering, such as piecewise affine controllers \cite{pasquale2020stabilizing}, can be readily considered by modifying the boundaries \eqref{eq_1_lb}-\eqref{eq_ub+} of the invariant set presented next.

By the conservation of flow, the dynamics of  the on-ramp queues $Q(t)$ and traffic densities $N(t)$ can be define below:
\begin{subequations}
\begin{align}
\dot{Q}^\mu_1(\phi)=&\alpha_1-r_1(q_1,n_1),\label{eq_G1} \\
\dot{N}^\mu_1(\phi)=&(r_1(q_1,n_1)-f_1^{\mu}(\phi))/l_1,\label{eq_H1} \\
\dot{Q}^\mu_k(\phi)=&\alpha_k-r_k^\mu(q_k,n_k), ~k\geq 2,\label{eq_G} \\
\dot{N}^\mu_k(\phi)=&(\beta_{k-1}f_{k-1}^\mu(\phi)+r_k^\mu(q_k,n_k)-f_k^\mu(\phi))/l_k, ~ k\geq2.\label{eq_H}
\end{align}
\label{eq_GH}%
\end{subequations}
For notional convenience, we define the vector fields $G^\mu(\phi):=[\dot{Q}_1^\mu,\cdots,\dot{Q}^\mu_K(\phi)]^{\mathrm{T}}$ and $H^\mu(\phi):=[\dot{N}_1^\mu,\cdots,\dot{N}^\mu_K(\phi)]^{\mathrm{T}}$.
One can show that $G^\mu$ and $H^\mu$ are continuous and bounded and that $Q(t)$ and $N(t)$ are continuous in time $t$ \cite{jin19tac}.

Finally, we present lower and upper boundaries of traffic densities used for control design. For a set $\mathcal{S}\times\mathcal{E}\subseteq\mathcal{S}\times\mathcal{Q}\times\mathcal{N}$, we say $n_k^*$ is a lower bound (resp. upper bound) of density $n_k$ over $\mathcal{S}\times\mathcal{E}$ if $H_k^{\mu}(\phi) \geq 0$ (resp. $H_k^{\mu}(\phi) \leq 0$) for any $\phi\in\{(s,q,n)\in\mathcal{S}\times\mathcal{E}|n_k = n_k^*\}$. We define lower bounds \eqref{eq_1_lb}-\eqref{eq_2_lb+} and upper bounds \eqref{eq_ub}-\eqref{eq_ub+} as follows:
\begin{subequations}
\begin{align}
    \underline n_1 :=& \frac{\min\{\alpha_1, F_1^{\max}\}}{v_1}, \label{eq_1_lb} \\
    \underline n_k :=& \frac{\min\{\beta_{k-1}\min\{v_{k-1}\underline{n}_{k-1}, F_{k-1}^{\min}\}+\alpha_k, F_k^{\max}\}}{v_k},k\geq 2, \label{eq_2_lb} \\
    {\uwave n}{_1} :=& \frac{\min\{U_1, F_1^{\max}\}}{v_1}, \label{eq_1_lb+} \\
    {\uwave n}{_k^\mu} :=& \min\{\nu_k^{\mu}, \frac{F_k^{\max}}{v_k}\},~k\geq2, \label{eq_2_lb+} \\
    \bar{n}_k :=&  n_k^{\mathrm{jam}} - \frac{F_k^{\min}}{w_k},~k\leq K, \label{eq_ub} \\
    \tilde{n}_{k}^\mu :=& n_k^{\mathrm{jam}} - \frac{1}{w_k} \min\{F_k^\min, \frac{R_{k+1}}{\beta_k}\}, k\leq K-1, \label{eq_ub+}
\end{align}
\end{subequations}
where $\nu_k^\mu$ is a solution of $$\beta_{k-1}\min\{v_{k-1}\underline{n}_{k-1
}, F_{k-1}^{\min}\}+r_k^{\mu}(1, \nu_k^\mu) = v_k \nu_k^\mu$$ 
and 
$$
R_{k+1}:=\min_{n_{k+1}\in[\underline{n}_{k+1},\tilde{n}_{k+1}^\mu]}
     w_{k+1}(n_{k+1}^{\mathrm{jam}} - n_{k+1}) - r_{k+1}^{\mu}(1, n_{k+1})$$
denotes the minimum receiving flow of cell $k+1$ for mainline traffic. Note that $\nu_k^\mu$ is unique because of the monotonicity of $\mu$. \eqref{eq_1_lb}-\eqref{eq_2_lb} denote the lower bound of $n_k$ if there is no queue in buffer $k$, while \eqref{eq_1_lb+}-\eqref{eq_2_lb+} mean the lower bound of $n_k$ if buffer $k$ has queuing flow. \eqref{eq_ub} represents the upper bound of $n_k$ if cell $k$ is free from downstream congestion, while \eqref{eq_ub+} indicates the upper bound of $n_k$ if the outflow from cell $k$ is affected by downstream congestion. The superscript $\mu$ indicates the dependency of boundaries on the control law $\mu$. 

\subsection{Problem formulation}
Our control objective is to stabilize and improve highway system performance in the face of capacity perturbations. For performance metric under control policy $\mu$, we consider \emph{time-averaged queue size} $\bar{Q}^\mu$ in the buffers (on-ramps and upstream mainline). This quantity is useful for our stability analysis. It is also a typical cost function of in control of stochastic queuing/fluid models \cite{dai95}, and is related to the practical performance metric \emph{vehicle hours traveled} \cite{papageorgiou1997alinea}.

We consider the following notion of stability: the buffer queues are \emph{stable} if there exists $Z<\infty$ such that for each initial condition $\phi(0)\in\mathcal S\times\mathcal Q\times\mathcal N$,
\begin{align}
    \limsup_{t\to\infty}\frac1t\int_{\tau=0}^t\mathrm E\Big[\sum_{k=1}^KQ_k^\mu(\tau)\Big]\mathrm{d}\tau
    \le Z,
    \label{eq_bounded}
\end{align}
where $Q_k^\mu(\tau)$ denotes the queue size in the $k$-th buffer under the control law $\mu$. Practically, this notion indicates that the time-averaged queue lengths in all buffers are bounded regardless of the initial condition. For stable buffer queues, the time-averaged queue size converges to $\bar{Q}^\mu$ almost surely (a.s.):
\begin{align}
    \lim_{t\to\infty}\frac1t\int_{\tau=0}^t\Big(\sum_{k=1}^KQ_k^\mu(\tau)\Big)\mathrm{d}\tau = \bar{Q}^\mu ~a.s.
    \label{eq_Qtau}
\end{align}
We now state our problem formulation:
\begin{align}
    (\mathrm{P}_0)~\min_{\mu} & ~ \bar{Q}^\mu\nonumber\\
    s.t.~& \text{every buffer queue is bounded on average}.
\end{align}
The problem $\mathrm{P}_0$ aims at improving the system performance subject to the stability constraint (11). Note that the objective function comprises queues both at on-ramps and upstream mainline. Shorter on-ramp queues indicates less queuing time at on-ramps whereas shorter upstream mainline queue implies less severe mainline congestion.

In the following, we apply the \emph{Foster-Lyapunov criterion} to express $\mathrm{P}_0$ mathematically. The criterion presents a sufficient condition for \eqref{eq_bounded} and is stated below:

\noindent{\bf Foster-Lyapunov criterion \cite{meyn93}}.
\emph{Consider a Markov process $X(t)$ with state space $\mathcal{X}$ and infinitesimal generator $\mathscr L$. If there exists a Lyapunov function $V:\mathcal{X}\to\mathbb R_{\ge0}$ and constants $c>0$, $d<\infty$ satisfying
\begin{align}
    \mathscr L V(x)\le-cg(x)+d,~\forall x\in\mathcal{X},
    \label{eq_drift}
\end{align}
then for each initial condition $x\in\mathcal X$
\begin{align}
    \limsup_{t\to\infty}\frac1t\int_{\tau=0}^t\mathrm E[g(X(t))]\mathrm{d}\tau\le \frac{d}{c}.
    \label{eq_d/c}
\end{align}
}

We refer to \eqref{eq_drift} as the \emph{drift condition}. Noting that the traffic model SS-CTM belongs to piecewise deterministic Markov process \cite{jin19tac}, we have
\begin{align}
    \mathscr L V(\phi) =& G^\mu(\phi)^{\mathrm{T}}\nabla_q V(\phi)
    +H^\mu(\phi)^{\mathrm{T}}\nabla_n V(\phi) \nonumber \\
    &+ \sum_{s'\in\mathcal S}\lambda_{s,s'}(V(s',q,n)-V(s,q,n))
    \label{eq_infini}
\end{align}
for any Lyapunov function $V:\mathcal{S}\times\mathcal{Q}\times\mathcal{N}\to\mathbb{R}_{\ge0}$. 

The criterion states that $g(X(t))$ is bounded on average with an upper boundary $d/c$ if the drift condition \eqref{eq_drift} holds for $g$. Practically, the choice of $g$ is application dependent. For our problem, we let 
\begin{equation}
 g(\phi)=\sum_{k=1}^Kq_k, \label{eq_g}
\end{equation}
i.e., the sum of buffer queue sizes. Hence, the Foster-Lyapunov criterion can be used to bound time-averaged buffer size.

Importantly, to ensure stability, we do not need to verify the drift condition \eqref{eq_drift} over $\mathcal{S}\times\mathcal{Q}\times\mathcal{N}$. We  only need to consider initial conditions in an \emph{invariant set} \cite{jin19tac}, denoted by $\mathcal M^\mu\subseteq\mathcal Q\times\mathcal N$, where the superscript $\mu$ indicates the dependency of $\mathcal M^\mu$ on $\mu$. An invariant set $\mathcal M^\mu$ is a closed set such that $(Q(t),N(t))\in\mathcal M^\mu$ for $t\geq0$ given any initial condition $\phi\in\mathcal S\times\mathcal M^\mu$.

Since it is hard to compute $\bar{Q}^\mu$ analytically, we choose the upper bound $d/c$ of $\bar{Q}^\mu$ as the control objective. Given a invariant set $\mathcal{M}^\mu$ and a Lyapunov function $V:\mathcal{S}\times\mathcal{M}^\mu\to\mathbb{R}_{\geq0}$, the problem $\mathrm{P}_0$ can be approximated as follows:
\begin{align}
    (\mathrm{P}_1)~&\min_{\mu;c,d>0}~\frac{d}{c}\nonumber\\
    & s.t.~ G^\mu(\phi)^{\mathrm{T}}\nabla_q V(\phi)
    +H^\mu(\phi)^{\mathrm{T}}\nabla_n V(\phi) \nonumber \\
    &\quad +\sum_{s'\in\mathcal S}\lambda_{s,s'}(V(s',q,n)-V(\phi))
    \leq -c \sum_{1\leq k\leq K} q_k + d,
    \nonumber \\
    &\quad\quad\quad\quad\quad\quad \forall\phi=(s,q,n)\in \mathcal{S}\times\mathcal M^\mu, \label{eq_P1_const}
\end{align}
where the constraint \eqref{eq_P1_const} is derived by plugging \eqref{eq_infini}-\eqref{eq_g} into \eqref{eq_drift}. The problem $\mathrm{P}_1$ models simultaneous control of multiple ramps (i.e., coordinated control design).

In general, $\mathrm{P}_1$ is suboptimal to the original problem $\mathrm{P}_0$, since the upper bound $d/c$ is a proxy for the mean queue length $\bar{Q}^\mu$. Still, directly solving $\mathrm{P}_1$ is challenging, because the objective function is non-convex and the constraint \eqref{eq_P1_const} can be non-convex for a general Lyapunov function $V(\phi)$.  Besides, the constraint \eqref{eq_P1_const} is infinite dimensional. To address these issues, we first consider localized ramp metering in Section III and reformulate $\mathrm{P}_1$ as a bi-level program that can be solved quickly. Then we consider coordinated ramp metering in Section IV and provide another bi-level program as well. For the large-scale coordinated control, we present a partially coordinated approach that decomposes the problem into a series of sub-problems. Our numerical example later demonstrates that this approach can provide a good quality solution.
\section{Localized control design}
\label{sec_distributed}

\begin{figure}[htb]
\centering
\includegraphics[width=0.4\linewidth]{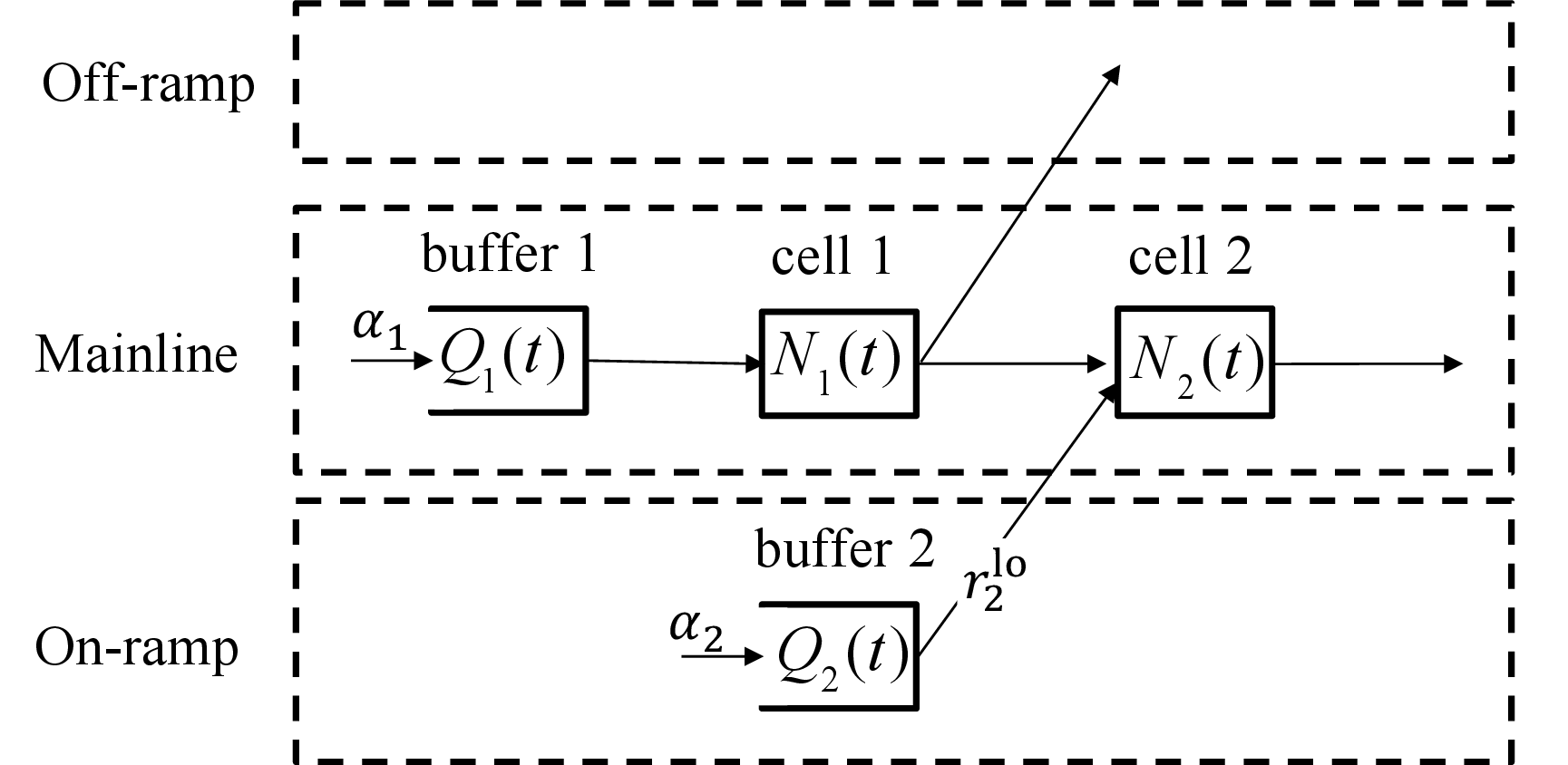}
\caption{Two-cell highway section.}
\label{fig_2cell}
\end{figure}

We begin by focusing on the design of localized ramp metering for a two-cell highway section illustrated by Fig.~\ref{fig_2cell}. Cells 1 and 2 (mainline) and buffer 2 (on-ramp) make up a typical merging area, which in practice may become a bottleneck due to merging conflicts and capacity perturbations \cite{asgharzadeh2020effect}. We assume that the outflow from cell 2 discharges with free-flow speed when cell 2 is congested. 
The control policy $\mu^{\mathrm{lo}}$ regulates the flow $r_2^{\mathrm{lo}}$ from buffer 2 to cell 2. Here we use the superscript ``$\mathrm{lo}$'' to indicate localized control. This ramp metering is localized because it focuses on control of the two-cell segment and is designed independently of its upstream or downstream ramp meters. In this section, we show that such a controller can be obtained by solving a bi-level program $\mathrm{P}_2$. To solve $\mathrm{P}_2$, we perform a grid search over the set of control parameters, and solve a quadratic program for each of candidate parameters.

Next, we present the stability condition of two-cell model and illustrate the localized control design.

\subsection{Main result}

Our approach consists of two steps: 1) finding a set of stabilizing controllers and 2) determining an ``optimal'' stabilizing controller. Essentially, we first convert the stability constraint \eqref{eq_P1_const} to an easy-to-check criterion for determining whether a localized controller is stabilizing. This is achieved by specifying the Lyapunov function $V(\phi)$ and the invariant set $\mathcal{M}^\mu$ in $\mathrm{P}_1$, both which are informed by the two-cell highway dynamics.

Intuitively, the highway section is stable if the time-averaged net flow during congestion period is negative. Here the net flow refers to inflow minus outflow, and the congestion indicates that there are queuing vehicles in buffers 1 or 2. Formally, given $\phi=(s,q,n)$, we define
\begin{subequations}
\begin{align}
    D_{1}^{{\mathrm{lo}}}(\phi)&:= \Big(G_1^{{\mathrm{lo}}}(\phi) +  l_1H_1^{{\mathrm{lo}}}(\phi)\Big) \nonumber \\
    & + \beta_1\Big(G_2^{{\mathrm{lo}}}(\phi) + \rho_2(n) l_2H_2^{{\mathrm{lo}}}(\phi)\Big), \label{eq_net_1} \\
    D_{2}^{{\mathrm{lo}}}(\phi)&:= \beta_1 \Big(G_1^{{\mathrm{lo}}}(\phi) +  l_1H_1^{{\mathrm{lo}}}(\phi)\Big) \nonumber \\
    & + \Big(G_2^{{\mathrm{lo}}}(\phi) + \rho_2(n) l_2H_2^{{\mathrm{lo}}}(\phi)\Big), \label{eq_net_2} 
\end{align}
\end{subequations}
where $\rho_2(n)$ is a weight function. We use the superscript ``$\mathrm{lo}$'' instead of $\mu$ to emphasize the dependence on localized control. We interpret $D_1^{{\mathrm{lo}}}(\phi)$ as a combination of weighted net flows of the cells downstream of buffer 1 (i.e., cells 1 and 2) and $D_2^{{\mathrm{lo}}}(\phi)$ as a weighted net flow of the cell downstream of buffer 2 (i.e., cell 2). To see this, let us consider $\rho_2(n)\equiv1$ and obtain 
\begin{align*}
    D_1^{\mathrm{lo}}(\phi)=&(1-\beta_1^2)(\alpha_1-f_1^{\mathrm{lo}}(\phi))+\beta_1(\beta_1\alpha_1+\alpha_2-f_2^{\mathrm{lo}}(\phi)), \\
    D_2^{\mathrm{lo}}(\phi) =&\beta_1\alpha_1+\alpha_2-f_2^{\mathrm{lo}}(\phi).
\end{align*}
Clearly, $\alpha_1-f_1^{\mathrm{lo}}(\phi)$ is the net flow of cell 1 and $\beta_1\alpha_1+\alpha_2-f_2^{\mathrm{lo}}(\phi)$ is the net flow of cell 2. Although the simple weight function leads to an intuitive definition of weighted net flows $D_1^{\mathrm{lo}}(\phi)$ and $D_2^{\mathrm{lo}}(\phi)$, these quantities neglect the on-ramp flow $r_2^{\mathrm{lo}}$ that merges into mainline traffic. Hence $\rho_2(n)\equiv1$ typically does not provide a tight stability result. Instead, we choose the following weight function:
\begin{equation}
    \rho_2(n) = (n_2-\underline{n}_2)/(\bar{n}_2-\underline{n}_2),
\end{equation}
where $\underline{n}_2$ and $\bar{n}_2$ are given by \eqref{eq_2_lb} and \eqref{eq_ub}, respectively. This new weight function ensures that $H_2^{{\mathrm{lo}}}(\phi)$ has more influence on the net flow when the traffic density $n_2$ is higher. Besides, it leads to an efficiently solvable control design problem.

Noting $\phi=(s,q,n)$, we let 
\begin{align}
    D_{k,s}^{{\mathrm{lo}}}(\mathcal{E}):=&\max_{(q,n)\in\mathcal{E}} ~D_k^{\mathrm{lo}}(\phi), ~k=1,2
    \label{eq_net_outflow}
\end{align}
denote the maximum net flow over a given traffic state set $\mathcal{E}$ under mode $s$ and the control policy $\mu^{\mathrm{lo}}$. Then we consider two sets of $(q_1,q_2,n_1,n_2)$:
\begin{subequations}
\begin{align}
   \mathcal{E}_1^{{\mathrm{lo}}} :=&(\{1\}\times\{0\}\times\{\uwave{n}{_1}\}\times[\underline{n}_2, \bar{n}_2]) \nonumber
    \\ & \cup(\{1\}\times\{1\}\times\{\uwave{n}{_1}\}\times[\uwave{n}{_2^{{\mathrm{lo}}}}, \bar{n}_2]), \label{eq_set_E1} \\
    \mathcal{E}_2^{{\mathrm{lo}}} :=&\{0\}\times\{1\}\times\{\underline n_1\} \times [\uwave{n}{_2^{{\mathrm{lo}}}}, \bar{n}_2], \label{eq_set_E2}
\end{align}
\end{subequations}
which is $\underline{n}_1$, $\uwave{n}{_1}$, $\underline{n}_2$, $\bar{n}_2$ and $\uwave{n}{_2^{{\mathrm{lo}}}}$ are given by \eqref{eq_1_lb}-\eqref{eq_ub+}. Note that $\mathcal{E}_k^{\mathrm{lo}}$, $k=1,2$, denotes a set of states with vehicles queuing in buffer $k$ (i.e., congestion in buffer $k$). It turns out that we only need to consider maximal net flows over these two sets.


We are now ready to state the main result of this section:
\begin{thm}
\label{thm_distributed}
Consider the two-cell highway model with an affine control policy $\mu^{\mathrm{lo}}$. Then, the highway model is stable if
\begin{equation}
    \bar D^{{\mathrm{lo}}} := \max_{k=1,2}\sum_{s\in\mathcal S}p_sD^{{\mathrm{lo}}}_{k,s}(\mathcal{E}_k^{{\mathrm{lo}}}) < 0,
    \label{eq_psDs}
\end{equation}
where $\{p_s;s\in\mathcal{S}\}$ is the solution to \eqref{eq_prob_a}-\eqref{eq_prob_c}.
Furthermore, if \eqref{eq_psDs} holds, the on-ramp queues are upper-bounded by
\begin{align}
    &\limsup_{t\to\infty}\frac1t\int_{\tau=0}^t\mathrm E[|Q_1(\tau)+Q_2(\tau)|]\mathrm{d}\tau
    \le -\frac{d}{ \bar{D}^{{\mathrm{lo}}}},
\end{align}
where $d\in(0,\infty)$ is a constant independent of $\mu^{\mathrm{lo}}$. \qed
\end{thm}

Theorem~\ref{thm_distributed} is proved by introducing the piecewise Lyapunov function $V^{\mathrm{lo}}(\phi)$ and the invariant set $\mathcal{M}^{\mathrm{lo}}$; see Appendix~B. Here ``piecewise'' accounts for finite modes of the Markov process. We refer to $\bar D^{{\mathrm{lo}}}$ as the \emph{mean drift}. In other words, the above result says that the two-cell system is stable and the buffer queue size is bounded on average if the mean drift is negative. The negative mean drift \eqref{eq_psDs} indicates that when there are queuing vehicles at buffer $k\in\{1,2\}$,  the expected net flow of its downstream cells is negative. In that case, the traffic densities at the downstream cells tend to decrease in the long term and do not always block the queue discharge from buffer $k$.


Now we present the controller design based on the stability criterion. Because $d$ is independent of $\mu^{\mathrm{lo}}$, minimizing the upper bound $-d/\bar{D}^{{\mathrm{lo}}}$ is equivalent to minimizing the negative mean drift $\bar{D}^{{\mathrm{lo}}}$. Thus, the objective function of $\mathrm{P}_1$ can be replaced by $\bar{D}^{{\mathrm{lo}}}$. 
The control design for the two-cell highway section can be formulated as
$$\min_{u_2,\kappa_2>0} \bar{D}^{\mathrm{lo}}~s.t.~ \bar{D}^{\mathrm{lo}} <0,$$
where $u_2$ and $\kappa_2$ are control parameters (see \eqref{eq_controller}). Substituting the traffic dynamics \eqref{eq_G1}-\eqref{eq_H} into \eqref{eq_psDs}, we can write the above problem as:
\begin{subequations}
\begin{align}
    (\mathrm{P}_2)~&\min_{ u_2,\kappa_2>0,\bar{D}<0} \bar{D}\nonumber \\
    s.t.~&\bar{D}\geq\sum_{s\in\mathcal{S}} p_s\max_{(q,n)\in\mathcal{E}_1^{{\mathrm{lo}}}}\alpha_1 + \alpha_2 - (1-\beta_1^2\rho_2(n)) f_1^{{\mathrm{lo}}}(\phi)
     \nonumber \\
    & - \beta_1(1-\rho_2(n))r_2^{\mathrm{lo}}(q_2,n_2)
    - \beta_1\rho_2(n) f_2^{{\mathrm{lo}}}(\phi), \label{eq_local_constr1} \\
    &\bar{D}\geq\sum_{s\in\mathcal{S}} p_s\max_{(q,n)\in\mathcal{E}_2^{{\mathrm{lo}}}}\beta_1\alpha_1+\alpha_2
    -\beta_1(1-\rho_2(n))f_1^{\mathrm{lo}}(\phi) \nonumber \\
    & -(1-\rho_2(n))r_2^{\mathrm{lo}}(q_2,n_2)-\rho_2(n)f_2^{\mathrm{lo}}(\phi), \label{eq_local_constr2} 
\end{align}
\end{subequations}
where $\bar{D}$ is an auxiliary variable. Clearly, if $\mathrm{P}_2$ is feasible we have $\bar{D}^*=\bar{D}^{\mathrm{lo}*}$ at optimality.

In the above bi-level optimization program, $u_2$, $\kappa_2$ and $\bar{D}$ are the upper-level decision variables and $q_1$, $q_2$, $n_1$ and $n_2$ are the lower-level decision variables. The upper-level variables are interlinked with the lower-level variables by multiplication in the non-convex lower-level objective function. Given $u_2$ and $\kappa_2$, the lower-level problem essentially involves a finite number of non-convex quadratic problems. Fortunately, modern solvers such as Gurobi \cite{gurobi} can quickly solve such problems to global optimality. Besides, in practice, $u_2$ and $\kappa_2$ take values from compact intervals. Thus, we can compute the designed controller via a grid search over that ranges of $u_2$ and $\kappa_2$.

The drift condition \eqref{eq_psDs} may not hold under high traffic demand, which means that we cannot conclude stability or the system might be unstable. In that case, the program $\mathrm{P}_2$ can become infeasible. In peak hours when demands exceed highway capacity, a reasonable strategy is to maximize the throughput \cite{coogan15}. Thus, we consider an alternative reformulation $\tilde{\mathrm{P}}_2$ that maximizes the mainline throughput:
\begin{subequations}
\begin{align}
    (\tilde{\mathrm{P}}_2)~&\max_{ u_2,\kappa_2,\tilde{\alpha}_1>0,\bar{D}<0} \tilde{\alpha}_1 \nonumber \\
    s.t.~&\bar{D}\geq\sum_{s\in\mathcal{S}} p_s\max_{(q,n)\in\mathcal{E}_1^{{\mathrm{lo}}}}\tilde{\alpha}_1 + \alpha_2 - (1-\beta_1^2\rho_2(n)) f_1^{{\mathrm{lo}}}(\phi)
     \nonumber \\
    & - \beta_1(1-\rho_2(n))r_2^{\mathrm{lo}}(q_2,n_2)
    - \beta_1\rho_2(n) f_2^{{\mathrm{lo}}}(\phi), \label{eq_local_constr1_max} \\
    &\bar{D}\geq\sum_{s\in\mathcal{S}} p_s\max_{(q,n)\in\mathcal{E}_2^{{\mathrm{lo}}}}\beta_1\tilde{\alpha}_1+\alpha_2
    -\beta_1(1-\rho_2(n))f_1^{\mathrm{lo}}(\phi) \nonumber \\
    & -(1-\rho_2(n))r_2^{\mathrm{lo}}(q_2,n_2)-\rho_2(n)f_2^{\mathrm{lo}}(\phi), \label{eq_local_constr2_max} 
\end{align}
\end{subequations}
where $\tilde{\alpha}_1$ denotes a lower bound of the mainline throughput by noting that Theorem 1 is a sufficient condition.

While Theorem~\ref{thm_distributed} is only a sufficient condition, it can yield a sufficient and necessary condition in special cases. Consider a particular case where the localized control decouples cells 1 and 2 in the following way:
\begin{equation}
    w_2(n_2^{\mathrm{jam}} - n_2) - r_2^{{\mathrm{lo}}}(1,n_2) \geq \beta_1F_1^{\max},~\forall n_2\in[\underline{n}_2,\bar{n}_2], \label{eq_decoupling}
\end{equation}
which means that the outflow from cell 1 is not impeded by the traffic in cell 2 and buffer 2. If the decoupled two-cell highway model satisfies the following conditions \eqref{eq_pr1_con1}-\eqref{eq_pr1_con3}, then our stability criterion is sufficient and necessary:
\begin{cor}
\label{prp_distributed_iff}
Consider the two-cell highway model such that
\begin{subequations}
\begin{align}
    U_1&\ge F_1^{\max}, \label{eq_pr1_con1} \\
    \beta_1F_1^{\min}+ r_2^{{\mathrm{lo}}}(0,n_2^c)
    &\ge F^\max_2,\label{eq_pr1_con2} \\
    \beta_1\alpha_1+ r_2^{{\mathrm{lo}}}(1,n_2^c)
    &\ge F^\max_2.\label{eq_pr1_con3}
\end{align}
\end{subequations}
Then, the highway model is stable if and only if
\begin{subequations}
\begin{align}
    &\alpha_1< \bar F_1, \label{eq_pr1_conc1} \\
    &\beta_1 \alpha_1 + \alpha_2< \bar F_2. \label{eq_pr1_conc2}
\end{align}
\end{subequations}
\qed
\end{cor}

Corollary~1 is proved in Appendix~E. The condition \eqref{eq_pr1_con1} states that the buffer 1 is not the bottleneck upstream of cell 1. The condition \eqref{eq_pr1_con2} indicates that cell 2 is the bottleneck if cell 1 is congested but there are no queues at the on-ramp. The condition \eqref{eq_pr1_con3} says that cell 2 is the bottleneck if cell 1 is uncongested but there is a queue at the on-ramp.

\subsection{Numerical example}

Consider the two-cell highway section in Fig.~\ref{fig_2cell} with the  parameters in Table~\ref{tab_para}. The capacity of cell 2 switches between 3000 veh/hr and 6000 veh/hr and the transition rates $\lambda_{1,2}$ and  $\lambda_{2,1}$ are set based on practical values (see the calibration in Table~\ref{tab_markovprocess}).

\begin{table}[hbt]
\centering
\footnotesize
\caption{Model parameters.}
\label{tab_para}
\begin{tabular}{@{}ccccc@{}}
\toprule
Parameter & Notation & Value & Unit \\ \midrule
Free-flow speed & $v_1$, $v_2$  & 100  & km/hr       \\
Congestion-wave speed  & $w_1$, $w_2$ & 25 & km/hr       \\
Jam density  & $n^{\mathrm{jam}}_1$, $n^{\mathrm{jam}}_2$  & 200, 300  & veh/km       \\ 
Mainline ratio & $\beta_1$ & 0.75 & - \\
Two-mode capacity & 
$\begin{bmatrix}
F_{1,1} & F_{2,1}\\
F_{1,2} & F_{2,2}
\end{bmatrix}$ & $\begin{bmatrix}
4000 & 6000 \\
4000 & 3000
\end{bmatrix}$  & veh/hr \\
Transition rate & $\lambda_{1,2}, \lambda_{2,1}$ & 0.9 & /hr \\
On-ramp capacity & $U_1$, $U_2$  & 4000, 1200  & veh/hr  \\
Demand  & $\alpha_1$, $\alpha_2$  & 3500, 600  & veh/hr  \\
\bottomrule
\end{tabular}
\end{table}

\begin{figure}[htbp]
\centering
\subfigure[$\bar{D}^{{\mathrm{lo}}}$ under various $(\kappa_2,u_2)$: the solution $(\kappa_2^*,u_2^*)$ to $\mathrm{P}_0$ and the solution $(\kappa_2^{\mathrm{lo}},u_2^{\mathrm{lo}})$ to $\mathrm{P}_2$.]{
\centering
\includegraphics[width=0.35\linewidth]{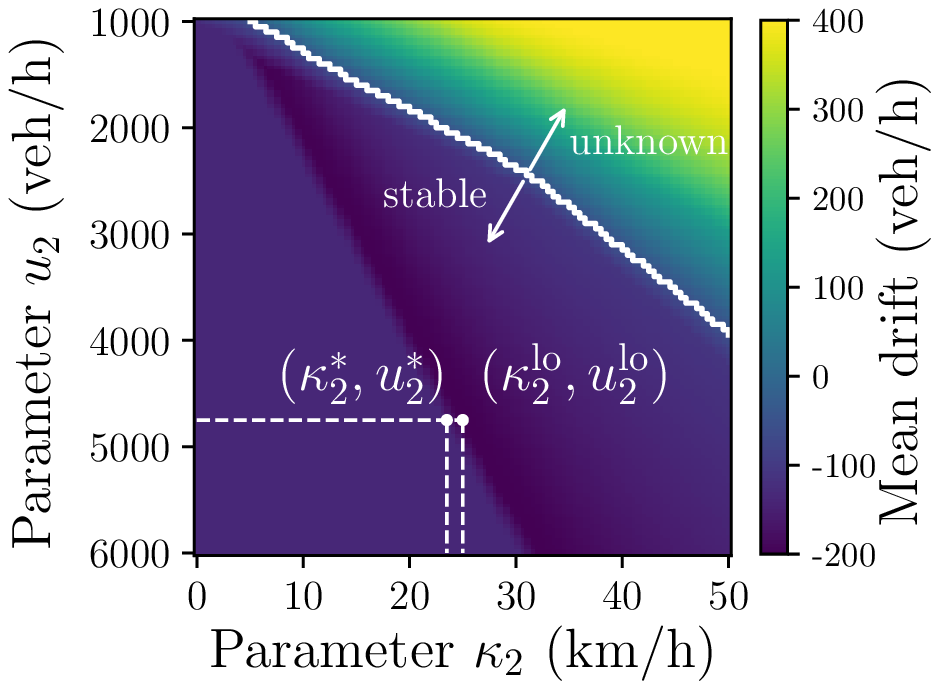}
\label{fig_drift_u_kappa}
}

\subfigure[$\bar{D}^{{\mathrm{lo}}}$ and $\hat{Q}^{{\mathrm{lo}}}$ under fixed $\kappa_2$.]{
\centering
\includegraphics[width=0.3\linewidth]{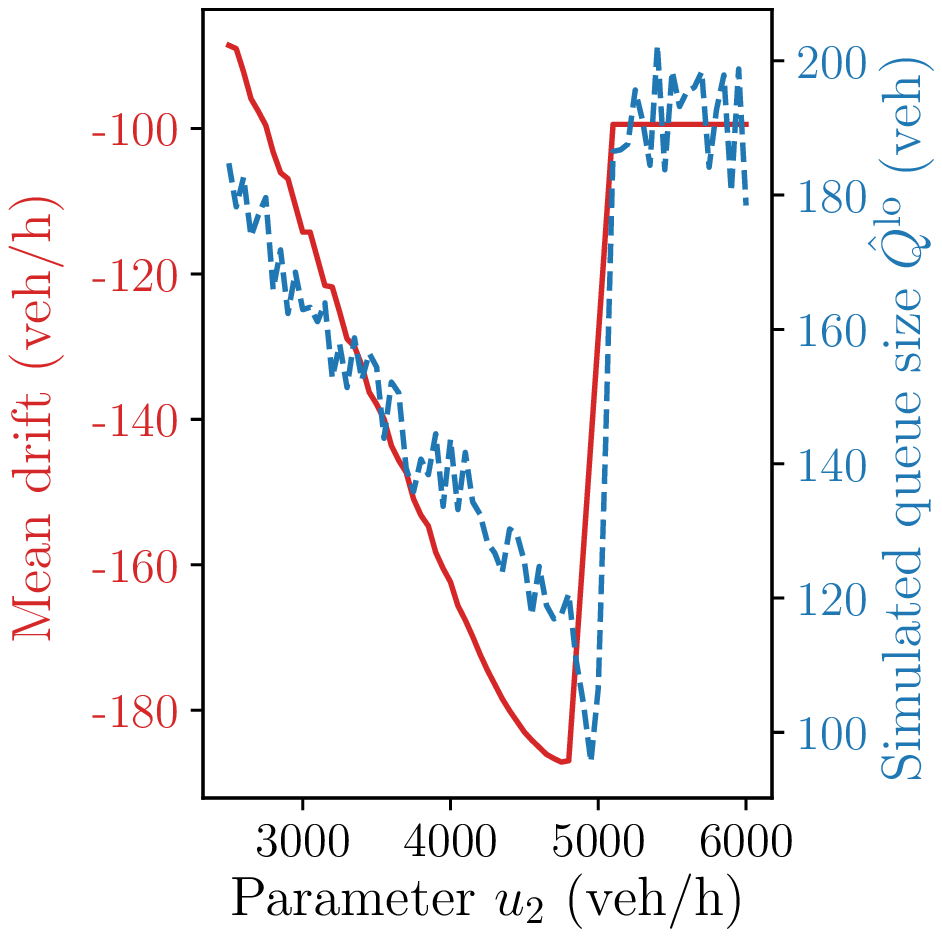}
\label{fig_drift_u}
}
\subfigure[$\bar{D}^{{\mathrm{lo}}}$ and $\hat{Q}^{{\mathrm{lo}}}$ under fixed $u_2$.]{
\centering
\includegraphics[width=0.3\linewidth]{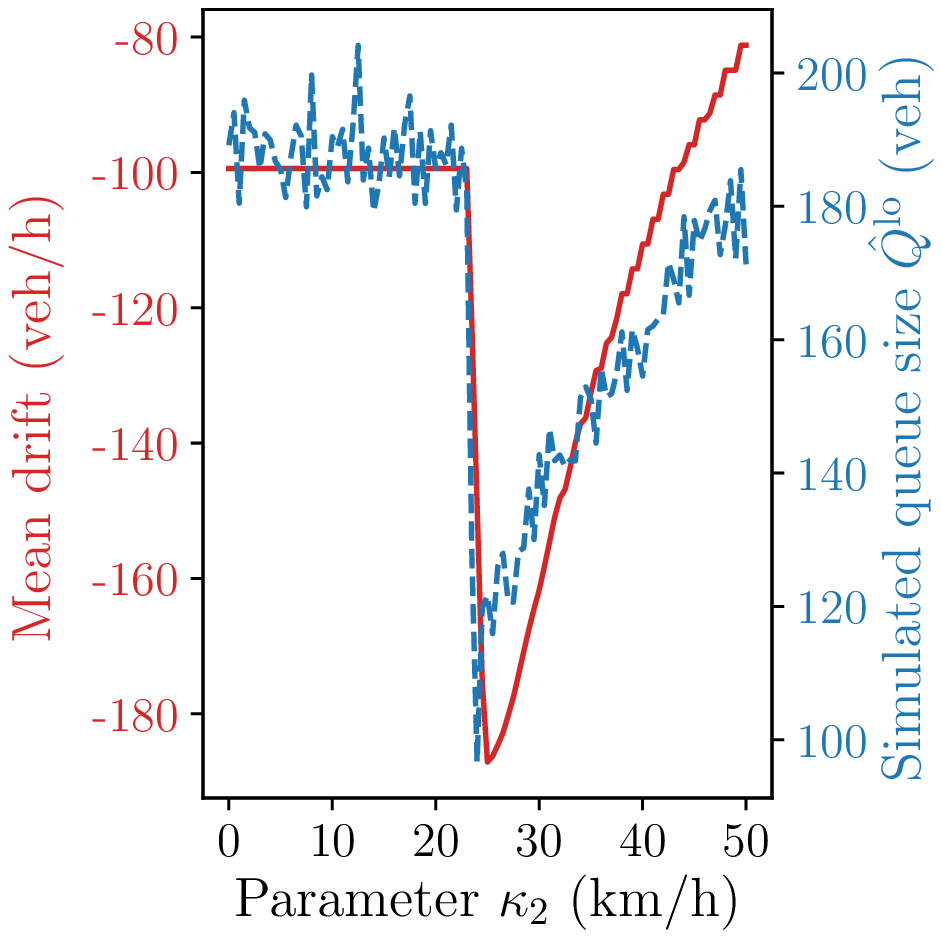}
\label{fig_drift_kappa}
}
\caption{Empirical relationship between mean drift $\bar{D}^{{\mathrm{lo}}}$ and simulated queue size $\hat{Q}^{{\mathrm{lo}}}$.}
\label{Dbar_v_Qbar}
\end{figure}

Fig.~\ref{fig_drift_u_kappa} presents the mean drift $\bar{D}^{{\mathrm{lo}}}$ for various values of $u_2$ and $\kappa_2$. We use Theorem~\ref{thm_distributed} to find the stable regime with negative $\bar{D}^{{\mathrm{lo}}}$. Since this theorem only states the sufficient condition for stability, it is unknown whether the regime with non-negative $\bar{D}^{{\mathrm{lo}}}$ is stable. But the result shows that $(\kappa_2^{\mathrm{lo}},u_2^{\mathrm{lo}})=(4750, 25)$ minimizes $\bar{D}^{{\mathrm{lo}}}$ and is the solution to $\mathrm{P}_2$. We also simulated the model numerically and obtained the average value of the total queue length $\hat{Q}^{\mu}=\frac{1}{T}\sum_{t=1}^T(q_1^{\mu}(t) + q_2^{\mu}(t))$, where $T$ denotes the number of simulation steps. We considered $T=10^6$ with a step size $\Delta_t=10$ seconds. The simulated value $\hat{Q}^{\mu}$ serves as an approximation of time-averaged value $\bar{Q}^{\mu}$ given by \eqref{eq_Qtau}. It was found that $(\kappa_2^*, u_2^*)=(4750, 24)$ minimizes  $\hat{Q}^{\mu}$ and can be seen as the solution to $\mathrm{P}_0$. Clearly, the example demonstrates that the program $\mathrm{P}_2$ can provide a good approximate solution for $\mathrm{P}_0$. 

We further illustrate the correspondence between $\bar{D}^{{\mathrm{lo}}}$ and $\hat{Q}^{{\mathrm{lo}}}$ by fixing one control parameter and tuning the other. Fig.~\ref{fig_drift_u} shows the case where $\kappa_2$ is fixed to 25 km/h and $u_2$ ranges between 2500 veh/h and 6000 veh/h, and Fig.~\ref{fig_drift_kappa} presents the other case where $u_2$ equals 4750 veh/h and $\kappa_2$ varies from 0 km/h to 50 km/h. As expected, we find that too small or large parameters lead to undesirable queue sizes. For example, aggressive ramp metering with small $u_2$ will increase the on-ramp queue size; lax ramp metering with large $u_2$ yields the control input $\mu_2^{\mathrm{lo}}(n_2)$ always larger than on-ramp flow $r_2^{\mathrm{lo}}(q_2,n_2)$ and thus does not exert any control effect on the highway section.

%

\begin{figure}[htbp]
\centering
\subfigure[Impacts on $u_2$.]{
\centering
\includegraphics[width=0.3\linewidth]{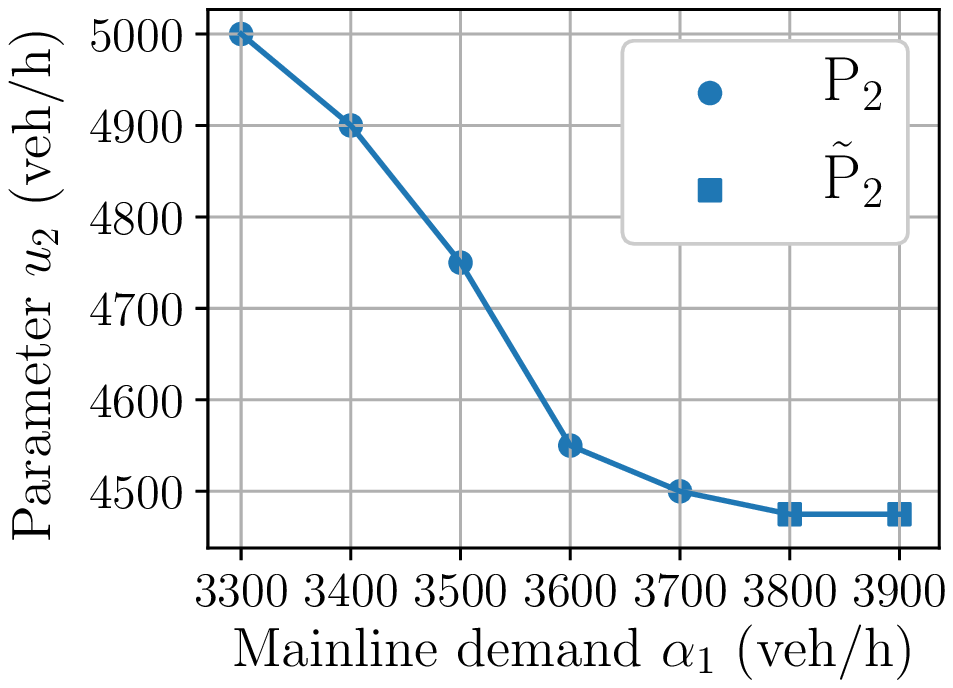}
\label{fig_twocell_u_analysis}
}
\subfigure[Impacts on $\kappa_2$.]{
\centering
\includegraphics[width=0.3\linewidth]{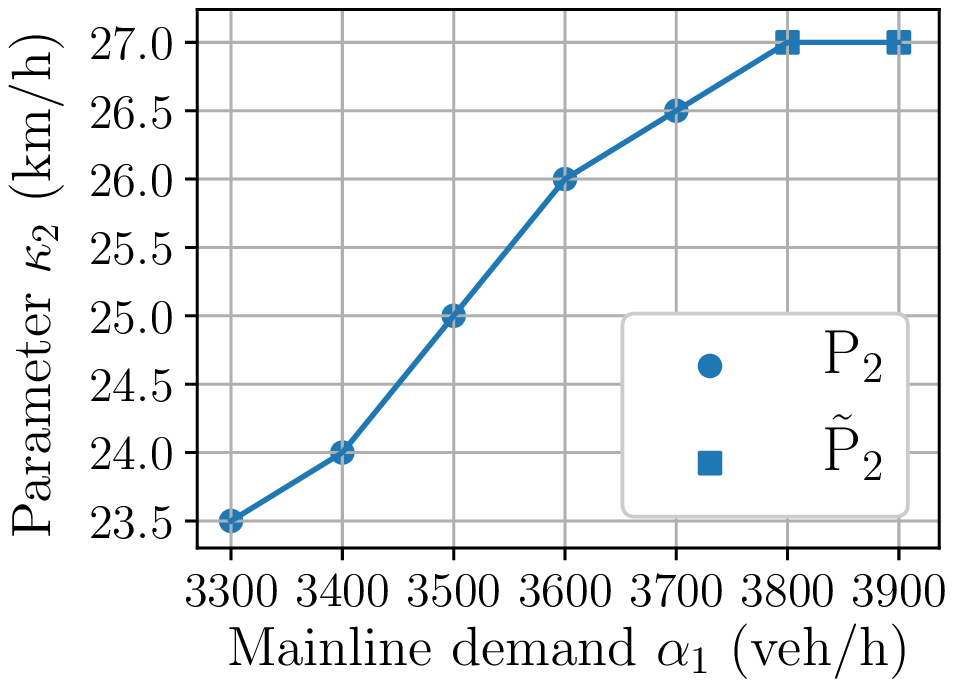}
\label{fig_twocell_kappa_analysis}
}
\caption{Impacts of demands on control parameters.}
\label{fig_twocell_analysis}
\end{figure}

We also explore the impact of demands on the control design. Fig.~\ref{fig_twocell_u_analysis} presents the relationship between the parameter $u_2$ and the mainline demand $\alpha_1$ given the on-ramp demand $\alpha_2 = 600$ veh/h and the parameter $\kappa_2=25$ km/h; Fig.~\ref{fig_twocell_kappa_analysis} presents the relationship between the parameter $\kappa_2$ and the mainline demand $\alpha_1$ given the on-ramp demand $\alpha_2 = 600$ veh/h and the parameter $u_2=4750$ veh/h. Note that the drift condition fails to hold given $\alpha_1\geq3800~\mathrm{veh/h}$. In that case, we turn to the program $\tilde{\mathrm{P}}_2$ to find the control parameters. We observe that $u_2$ tends to decrease and $\kappa_2$ tends to increase when $\alpha_1$ increases. This finding indicates that our control design yields sensible controllers that further restrict on-ramp flows when the mainline congestion worsens.

\begin{figure}[htbp]
    \centering
    \subfigure[Impacts of mainline demands.]{
    \centering
    \includegraphics[width=0.3\linewidth]{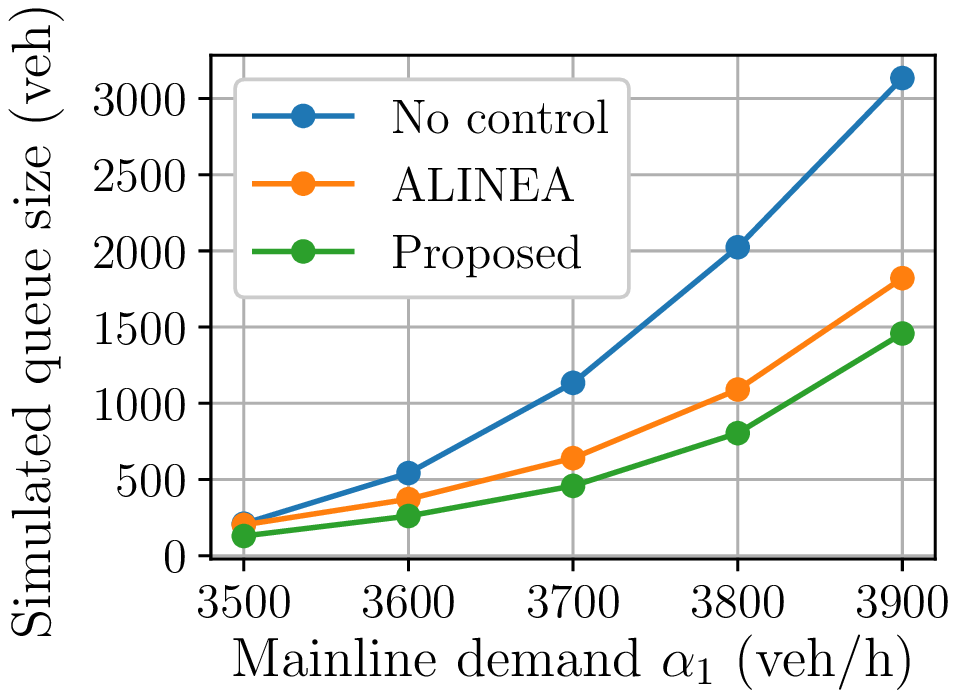}
    \label{fig_twocell_comp_demand}
    }
    \subfigure[Impacts of capacity fluctuation.]{
    \centering
    \includegraphics[width=0.3\linewidth]{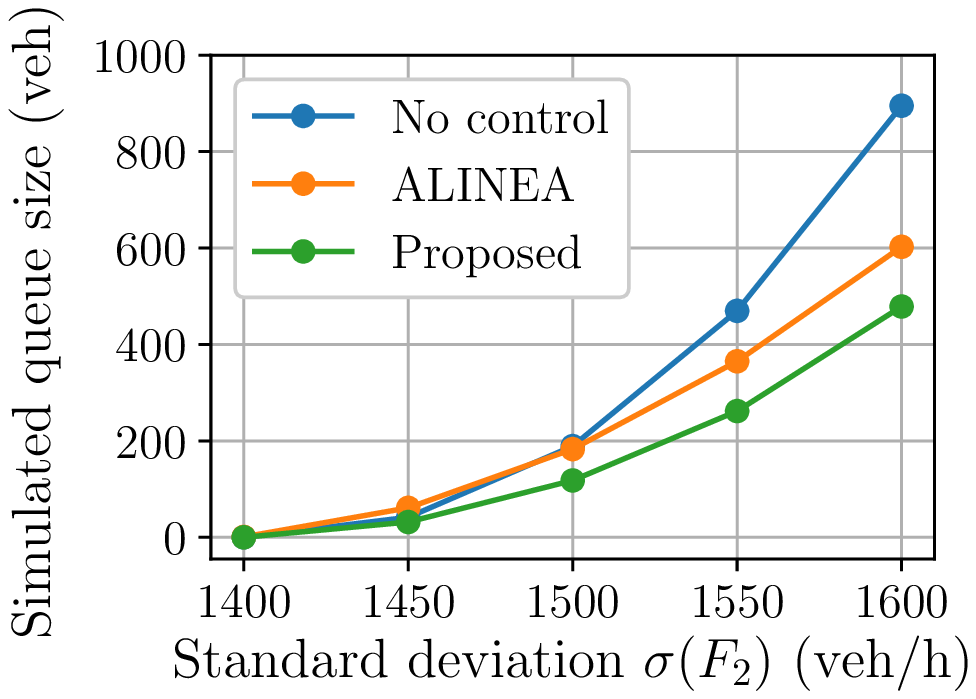}
    \label{fig_twocell_comp_std}
    }
    \label{fig_twocell_comp}
    \caption{Control comparison given different mainline demand and perturbed capacity.}
\end{figure}

Next, we compare the proposed method with the classical control ALINEA that regulates on-ramp flows $\mu_k(t)$ via the following feedback:
\begin{equation}
\mu_k(t) = \mu_k(t-1) - K_R(n_k(t) -n_k^c),
\end{equation}
where $K_R$ is the control gain, $n_k(t)$ is the traffic density of cell $k$ at time step $t$ and $n_k^c$ is the nominal critical density of cell $k$. ALINEA is essentially an integral controller and the literature suggested that the gain $K_R=40$ (km$\cdot$lane)/h is a reasonable choice \cite{wang2014local}. Fig.~\ref{fig_twocell_comp_demand} illustrates the comparison given different mainline demands $\alpha_1$. Both ALINEA and our method contribute to the marginal reduction of queue sizes when $\alpha_1$ is low. However, the ramp metering begins to work as the increasing mainline demand leads to long buffer queues; particularly our controller outperforms ALINEA in terms of shortening the queue size. Fig.~\ref{fig_twocell_comp_std} shows another comparison under perturbations in the mainline capacity $F_{2,2}$. Here we use $\sigma(F_2)$, standard deviation of $[F_{2,1}, F_{2,2}]$, as a proxy for capacity fluctuation. We again obtain similar performance improvement as Fig.~\ref{fig_twocell_comp_demand}. These numerical results suggest the benefit of our method in comparison to ALINEA in the face of stochastic capacity fluctuations. 

\section{Coordinated control design}
\label{sec_cordinated}

In principle, one can apply the localized control design to independently meter inflows through each on-ramp. However, this approach ignores the influence of downstream congestion. In practice, traffic jams can easily spill back and block the upstream on-ramps, especially when on-ramps are located closely. To address this issue, we propose a partially coordinated control design that sequentially determines ramp meters from the downstream to the upstream. Our main result of this section shows that the partial coordination can be addressed by sequentially solving optimization problems $\mathrm{P}_4^k$ for $k=K,K-1,\cdots,2$, where each $\mathrm{P}_4^k$ can be solved analogously to $\mathrm{P}_2$.

We first discuss the fully coordinated control design problem that simultaneously designs multiple ramp controllers, because it motivates the partially coordinated design problem. Next we illustrate the latter's benefits via a numerical example.

\subsection{Main result}
We consider the fully coordinated ramp metering $\mu^{\mathrm{fc}}$ applied to the $K$-cell highway shown in Fig.~\ref{fig_ctm}. Here the superscript ``$\mathrm{fc}$'' indicates full coordination. To state the results, we first generalize the weighted net flows \eqref{eq_net_1}-\eqref{eq_net_2} for the $K$-cell highway as follows:
\begin{align}
    D_k^{\mathrm{fc}}(\phi) :=& \sum_{j=1}^{k-1}\gamma_{j,k}\Big(G_j^{\mathrm{fc}}(\phi) + l_j\rho_j(n) H_j^{\mathrm{fc}}(\phi)\Big) \nonumber \\
    &+\sum_{j=k}^{K}\gamma_{k,j}\Big(G_j^{\mathrm{fc}}(\phi)+l_j\rho_j(n) H_j^{\mathrm{fc}}(\phi)\Big),
\end{align}
where
\begin{equation}
    \gamma_{j,k}:=
    \begin{cases}
     \prod_{\ell=j}^{k-1}\beta_{\ell} & j<k, \\
     1 & j=k,
    \end{cases}
     \label{eq_gamma} 
\end{equation}
and $\rho_j(n)$ is the weight function. We can interpret $D_k^{\mathrm{fc}}(\phi)$ as a combination of weighted net flows (inflow minus outflow) of cells downstream of buffer $k$, namely cells $k,k+1,\cdots,K$. Note that $\gamma_{j,k}$ for $j<k$ denotes the proportion of the mainline inflow of cell $k$ to the outflow of cell $j$. Inspired by the discussion in Section III, we generally consider $\rho_1\equiv1$ and $$\rho_j(n):=(n_j - \underline{n}_j)/(\tilde{n}_j^{\mathrm{fc}} - \underline{n}_j),~j=2,\cdots,K.$$ Recalling our notation of hybrid state $\phi=(s,q,n)$, we denote the corresponding maximum net flow over a given traffic state set $\mathcal{E}$ by
\begin{equation}
    D_{k,s}^{\mathrm{fc}}(\mathcal{E}):=
\max_{(q,n)\in\mathcal{E}} D_k^{\mathrm{fc}}(\phi), ~k=1,\cdots,K.
\end{equation}
We denote by
\begin{equation}
    \mathcal{E}^{{\mathrm{fc}}}_k:=\bigcup_{
    \substack{ i_k=1 \\
    i_j\in\{0,1\},1\leq j\leq K, j\neq k }} \prod_{j=1}^K
    \mathcal{I}_{j,i_j}^{{\mathrm{fc}}}\times
    \prod_{j=1}^K  \mathcal{N}_{j,i_j}^{{\mathrm{fc}}}
\end{equation}
a set of the state $(q_1,\cdots,q_K,n_1,\cdots,n_K)$ with queuing traffic in buffer $k$, where 
\begin{align*}
    \mathcal{I}_{j,i_j}^{{\mathrm{fc}}}:= \begin{cases}
    \{0\} & i_j=0 \\
    \{1\} & i_j=1
    \end{cases},~
    \mathcal{N}_{j,i_j}^{{\mathrm{fc}}}:= \begin{cases}
    [\underline{n}_j, \tilde{n}{_j^{{\mathrm{fc}}}}] & i_j=0 \\
    [{\uwave n}{_j^{{\mathrm{fc}}}}, \tilde{n}{_j^{{\mathrm{fc}}}}] & i_j=1
    \end{cases}.
\end{align*}
Here $i_j$ denotes whether there is a queue in buffer $j$.



We can now state the main result of this section:
\begin{thm}
\label{thm_centralized}
Consider a $K$-cell highway with an affine control policy $\mu^{\mathrm{fc}}$. Then, the highway model is stable if
\begin{align}
    \bar D^{{\mathrm{fc}}}:=\max_{k=1,\cdots,K}\sum_{s\in\mathcal S}p_s D_{k,s}^{{\mathrm{fc}}}(\mathcal{E}_k^{{\mathrm{fc}}})<0,
    \label{eq_psDs2}
\end{align}
where $\{p_s;s\in\mathcal{S}\}$ is the solution to \eqref{eq_prob_a}-\eqref{eq_prob_c}. Furthermore, if \eqref{eq_psDs2} holds, the buffer queues are upper-bounded by
\begin{align}
    &\limsup_{t\to\infty}\frac1t\int_{\tau=0}^t\mathrm E\Big[\sum_{k=1}^K Q_k(\tau)\Big]\mathrm{d}\tau
    \le -\frac{d}{\bar D^{{\mathrm{fc}}}}, \label{eq_ub_co}
\end{align}
where $d\in(0,\infty)$ is a constant independent of $\mu^{\mathrm{fc}}$. \qed
\end{thm}

The proof of Theorem~\ref{thm_centralized} is based on the piecewise Lyapunov function $V^{\mathrm{fc}}(\phi)$ and the invariant set $\mathcal{M}^{\mathrm{fc}}$; see Appendix~C. Recall from \eqref{eq_controller} that $\mu^{\mathrm{fc}}$ is parameterized by $u,\kappa\in\mathbb{R}^{K-1}_{>0}$; thus we can apply Theorem~\ref{thm_centralized} to recast $\mathrm{P}_1$ as $\mathrm{P}_3$:
\begin{align}
    (\mathrm{P}_3)&\min_{u,\kappa\in\mathbb{R}_{>0}^{K-1},\bar{D}<0}\bar{D} \nonumber \\
     s.t.~&\bar{D}\geq\sum_{s\in\mathcal{S}} p_s\max_{(q,n)\in\mathcal{E}_k^{{\mathrm{fc}}}} \bigg(
    \sum_{j=1}^{K}\tilde{\gamma}_{j,k}\alpha_j \nonumber \\ 
    & - \sum_{j=1}^{K}\Big(\tilde{\gamma}_{j,k}\rho_j(n)-\beta_{j}\tilde{\gamma}_{j+1,k}\rho_{j+1}(n)\Big)f_j^{\mathrm{fc}}(\phi) \nonumber \\
    & - \sum_{j=1}^K \tilde{\gamma}_{j,k}\Big(1-\rho_j(n)\Big) r_j^{{\mathrm{fc}}}(q_j,n_j)\bigg), ~1\leq k\leq K, \label{eq_p3_constr}
\end{align}
where $\bar{D}$ is an auxiliary variable analogous to that in $\mathrm{P}_2$, $\rho_{K+1}(n):=0$, $\tilde{\gamma}_{j,k}:=\gamma_{j,k}$ for $j\leq k$ and $\tilde{\gamma}_{j,k}:=\gamma_{k,j}$ for $j>k$. The constraints \eqref{eq_p3_constr} are derived by plugging \eqref{eq_G1}-\eqref{eq_H} into \eqref{eq_psDs2}. If $\mathrm{P}_3$ is infeasible, we can reformulate it as a throughput-maximization problem analogous to $\tilde{\mathrm{P}}_2$. The fully coordinated control design problem is also a bi-level program with non-convex lower-level optimization. Solving it is more challenging than $\mathrm{P}_2$, since both the lower- and upper-level decision variables increase with the scale of highway stretch. Consequently, it is even harder to find the globally optimal $u^{{\mathrm{fc}}}$ and $\kappa^{{\mathrm{fc}}}$ for realistic instances.

We propose a partially coordinated approach to address the computational issue of solving $\mathrm{P}_3$. 
The superscript ``$\mathrm{pc}$'' indicates partial coordination. The design of the control law $\mu^{\mathrm{pc}}$ follows the procedure: we first design the controller $\mu^{\mathrm{pc}}_K$ for ramp $K$, and then design the controller $\mu^{\mathrm{pc}}_{K-1}$ for ramp $K-1$ based on the controller $\mu^{\mathrm{pc}}_K$, and so forth until all $K-1$ controllers are obtained. 

To formulate the problem mathematically, we define
\begin{subequations}
\begin{align}
    D_k^{\mathrm{pc}}(\phi) := & \sum_{j=1}^{k-1}\gamma_{j,k}\Big(G_j^{\mathrm{pc}}(\phi) + l_j H_j^{\mathrm{pc}}(\phi)\Big) \nonumber \\
    &+\sum_{j=k}^{K}\gamma_{k,j}\Big(G_j^{\mathrm{pc}}(\phi) + l_j\rho_j(n) H_j^{\mathrm{pc}}(\phi)\Big), \label{eq_net_pc} \\
    D_{k,s}^{{\mathrm{pc}}}(\mathcal{E}):=&\max_{(q,n)\in\mathcal{E}}D_k^{\mathrm{pc}}(\phi), \label{eq_Dmax_pc} 
\end{align}
\end{subequations}
where $\phi=(s,q,n)$ and $\mathcal{E}$ is a given traffic state set. $D_k^{\mathrm{pc}}(\phi)$ (resp. $D_{k,s}^{{\mathrm{pc}}}(\mathcal{E})$) has a similar meaning to that of $D_k^{\mathrm{fc}}(\phi)$ (resp. $D_{k,s}^{{\mathrm{fc}}}(\mathcal{E})$). Note that we let
\begin{equation*}
    \rho_j(n) = \begin{cases}
     1 & j=1,2,\cdots,k-1, \\
     (n_j-\bar{n}_j)/(\tilde{n}_j^{\mathrm{pc}}-\bar{n}_j) & j=k,k+1,\cdots,K
    \end{cases}
\end{equation*}
in $D_k^{\mathrm{pc}}(\phi)$, which is different from $D_k^{\mathrm{fc}}(\phi)$. It implies that when computing the net flow, the partially coordinated control design treats the upstream highway along with on-ramps as a single cell. We also define a set of $(q_k,q_{k+1},\cdots,q_K,n_{k-1}, n_k,\cdots,n_K)$ for $k=2,\cdots,K$:
\begin{align}
    \mathcal{E}_k^{{\mathrm{pc}}}
    =& \bigcup_{i_{k+1},\cdots,i_K\in\{0,1\}}\{1\}\times\prod_{j=k+1}^K \mathcal{I}_{j,i_{j}}^{{\mathrm{pc}}} \nonumber \\
    &\times \{\underline{n}_{k-1}\}\times[\uwave{n}{_k^{{\mathrm{pc}}}},\tilde{n}_k^{{\mathrm{pc}}}]\times\prod_{j=k+1}^K\mathcal{N}_{j,i_{j}}^{{\mathrm{pc}}},
    \label{eq_E_pc}
\end{align}
where $\mathcal{I}_{j,i_j}^{{\mathrm{pc}}}$ and $\mathcal{N}_{j,i_j}^{{\mathrm{pc}}}$ are defined in a similar way to $\mathcal{I}_{j,i_j}^{{\mathrm{fc}}}$ and $\mathcal{N}_{j,i_j}^{{\mathrm{fc}}}$. The definition above indicates that when designing the controller $\mu^{\mathrm{pc}}_k$, the partial coordination only considers the cell $k-1$ among its upstream cells and buffers.

The controller $\mu^{\mathrm{pc}}_k$ parameterized by $u_k$ and $\kappa_k$ can be obtained by solving
\begin{align}
    (\mathrm{P}_4^k)~&\min_{u_k,\kappa_k>0, \bar{D}_k<0} \bar{D}_k \nonumber \\
    s.t.~&\bar{D}_{k} \geq \sum_{s\in\mathcal{S}}p_s\max_{(q,n)\in\mathcal{E}_k^{{\mathrm{pc}}}} \Big(\sum_{j=1}^{K}\tilde{\gamma}_{j,k} \alpha_j \nonumber \\
    & - \sum_{j=k-1}^{K}\Big(\tilde{\gamma}_{j,k}\rho_j(n)-\beta_j\tilde{\gamma}_{j+1,k}\rho_{j+1}(n)\Big)f_j^{\mathrm{pc}}(\phi) \nonumber \\
    & - \sum_{j=k-1}^K \tilde{\gamma}_{j,k}\Big(1-\rho_j(n)\Big) r_j^{{\mathrm{pc}}}(q_j,n_j)\bigg),
    \label{eq_constr_P4k}
\end{align}
where $\bar{D}_k$ is an auxiliary variable similar to that in $\mathrm{P}_2$.
For a particular $\mathrm{P}_4^k$, the boundaries $\tilde{n}_{k}^{{\mathrm{pc}}},\cdots, \tilde{n}_{K-1}^{{\mathrm{pc}}}$ and $\uwave{n}{_{k+1}^{{\mathrm{pc}}}},\cdots,\uwave{n}{_{K}^{{\mathrm{pc}}}}$ are determined by the controllers $\mu_{k+1}^{\mathrm{pc}},\cdots,\mu_{K}^{\mathrm{pc}}$ that have been already computed while $\uwave{n}{_{k}^{{\mathrm{pc}}}}$ relies on the controller $\mu_{k}^{\mathrm{pc}}$ to be solved. Importantly, the dimension of the upper-level decision variables $u_k,\kappa_k,\bar{D}_k$ in $\mathrm{P}_4^k$ does not grow with the scale of highway stretch; thus, these problems can be sequentially solved analogous to $\mathrm{P}_2$.

We also study the stability condition of the partially coordinated control design since $\mathrm{P}_4^k$, $k=2,\cdots,K$, is based on simplifications of $\mathrm{P}_3$. The stability condition is stated as follows:
\begin{thm}
\label{thm_cascaded}
Consider a $K$-cell highway section with an affine control policy $\mu^{\mathrm{pc}}$. Then, the highway model is stable if
\begin{align}
    \sum_{s\in\mathcal S}p_s D_{k,s}^{{\mathrm{pc}}}(\mathcal{E}_k^{{\mathrm{pc}}})<0,~ k=1,2,\ldots,K,
    \label{eq_psDsk}
\end{align}
where $\{p_s;s\in\mathcal{S}\}$ is the solution to \eqref{eq_prob_a}-\eqref{eq_prob_c}, and $D_{k,s}^{{\mathrm{pc}}}(\mathcal{E}_k^{{\mathrm{pc}}})$ is given by \eqref{eq_Dmax_pc} and \eqref{eq_E_pc}. \qed
\end{thm}

The proof of Theorem~\ref{thm_cascaded} is given in Appendix~D. Note that if $\mathrm{P}_4^k$, $k=2,\cdots,K$, are all feasible, the verification of stability condition only requires to check \eqref{eq_psDsk} for $k=1$, where $\mathcal{E}_1^{\mathrm{pc}}$ is given by
\begin{align*}
    \mathcal{E}_1^{{\mathrm{pc}}}
    = \bigcup_{i_{2},\cdots,i_K\in\{0,1\}}\{1\}\times\prod_{j=2}^K \mathcal{Q}_{j,i_{j}}^{{\mathrm{pc}}} \times[\uwave{n}{_1^{{\mathrm{pc}}}},\tilde{n}_1^{{\mathrm{pc}}}]\times\prod_{j=2}^K\mathcal{N}_{j,i_{j}}^{{\mathrm{pc}}}.
\end{align*}
If $\mathrm{P}_4^k$ is infeasible for some $k$ or \eqref{eq_psDsk} fails to holds, we can consider the partially  coordinated control that maximizes the throughput, which is analogous to $\mathrm{\tilde{P}}_2$.

So far, we have discussed the localized, fully coordinated and partially coordinated control designs. The following case identifies the conditions when all three designs are equivalent.
\begin{cor} \label{prp_equi}
Consider a $K$-cell highway that is decoupled by a control policy $\mu^{\mathrm{pc}}$ such that
\begin{align}
    &w_{k+1}(n^{\mathrm{jam}}_{k+1} - n_{k+1}) - r_{k+1}^{{\mathrm{pc}}}(1, n_{k+1}) \geq \beta_k F_k^{\max} \label{eq_thm3_con} \\
    &\quad\quad \forall n_{k+1}\in [\underline{n}_{k+1},\bar{n}_{k+1}], ~k=1,\cdots,K-1.  \nonumber
\end{align}
If the controlled on-ramp flows satisfy
\begin{subequations}
\begin{align}
    U_1&\ge F_1^{\max}, \label{eq_pr2_con1} \\
    \beta_{k-1} F_{k-1}^{\min}+ r_k^{{\mathrm{pc}}}(0,n_k^c)
    &\ge F^\max_k,\label{eq_pr2_con2} \\
    \sum_{i=1}^{k-1}\gamma_{i,k}\alpha_i+ r_k^{{\mathrm{pc}}}(1,n_k^c)
    &\ge F^\max_k, \label{eq_pr2_con3}
\end{align}
\end{subequations}
then $\mu^{\mathrm{pc}}$ parameterized by $u^{\mathrm{pc}},\kappa^{\mathrm{pc}}\in\mathbb{R}_{>0}^{K-1}$ is also an optimal solution to $\mathrm{P}_2$ and $\mathrm{P}_3$ with $\rho_j\equiv1$ for $j=1,2,\cdots,K$. \qed
\end{cor}

The proof of Corollary~2 is given in Appendix~F. The decoupling condition \eqref{eq_thm3_con} has the same interpretation as \eqref{eq_decoupling}, and the conditions \eqref{eq_pr2_con1}-\eqref{eq_pr2_con3} are also analogous to \eqref{eq_pr1_con1}-\eqref{eq_pr1_con3}. By saying that $\mu^{\mathrm{pc}}$ is optimal to $\mathrm{P}_2$, we mean that $u_k^{\mathrm{pc}}$ and $\kappa_k^{\mathrm{pc}}$ are optimal if we apply $\mathrm{P}_2$ to the two-cell section made up of cells $k-1$ and $k$. 
In general, the decoupling condition does not hold. Thus one can expect the coordinated control to be superior to the localized control. The extent to which full coordination outperforms the partial coordination depends on the problem instance, as explained next.

\subsection{Numerical example}
\begin{figure}[htb]
\centering
\includegraphics[width=0.4\linewidth]{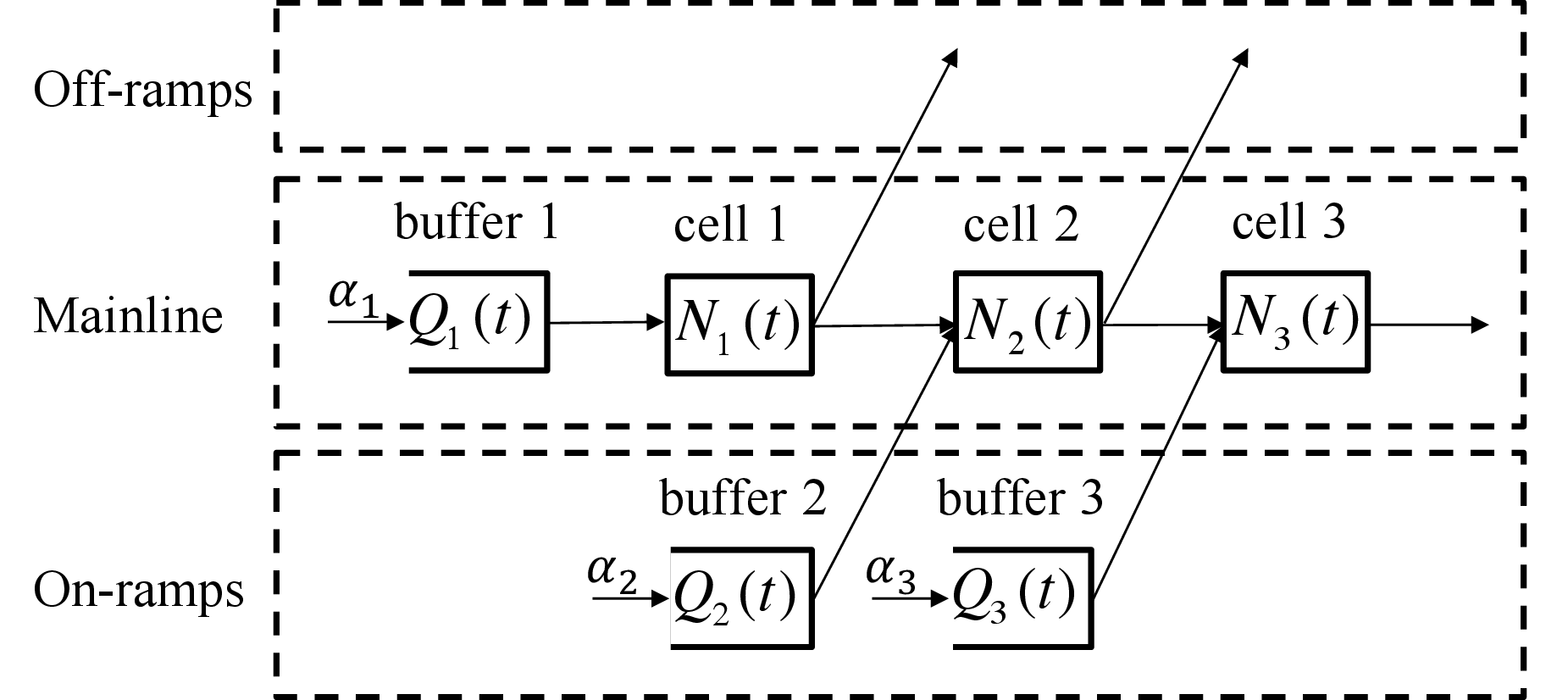}
\caption{Three-cell highway section.}
\label{fig_3cell}
\end{figure}

We now consider a three-cell highway section (see Fig.~\ref{fig_3cell}) with the  parameters in Table~\ref{tab_para_3}. In addition to cell 2, cell 3 also suffers capacity fluctuation between 3000 veh/hr and 6000 veh/hr. For the convenience of analysis, we let both $\kappa_2$ and $\kappa_3$ equal 25 km/h and optimize $u_2$ and $u_3$.

\begin{table}[hbtp]
\centering
\scriptsize
\caption{Model parameters.}
\label{tab_para_3}
\begin{tabular}{@{}ccccc@{}}
\toprule
Parameter & Notation & Value & Unit \\ \midrule
Free-flow speed & $v_1$, $v_2$, $v_3$  & 100  & km/hr       \\
Congestion-wave speed  & $w_1$, $w_2$, $w_3$ & 25 & km/hr       \\
Jam density  & $n^{\mathrm{jam}}_1$, $n^{\mathrm{jam}}_2$, $n^{\mathrm{jam}}_3$ & 200, 300, 300  & veh/km       \\ 
Mainline ratio & $\beta_1$, $\beta_2$ & 0.75, 0.6 & - \\
Two-mode capacity & 
$\begin{bmatrix}
F_{1,1} & F_{1,2} \\
F_{2,1} & F_{2,2} \\
F_{3,1} & F_{3,2}
\end{bmatrix}$ & $\begin{bmatrix}
4000 & 4000 \\
6000 & 3000 \\
6000 & 2500
\end{bmatrix}$  & veh/hr \\
Transition rate & $\lambda_{1,2}, \lambda_{2,1}$ & 0.9 & /hr \\
On-ramp capacity & $U_1$, $U_2$, $U_3$  & 4000, 1200, 1200  & veh/hr  \\
Demand  & $\alpha_1$, $\alpha_2$, $\alpha_3$  & 3500, 600, 800  & veh/hr  \\
\bottomrule
\end{tabular}
\end{table}

\begin{figure}[htbp]
\centering
\subfigure[Mean drift.]{
\centering
\includegraphics[width=0.3\linewidth]{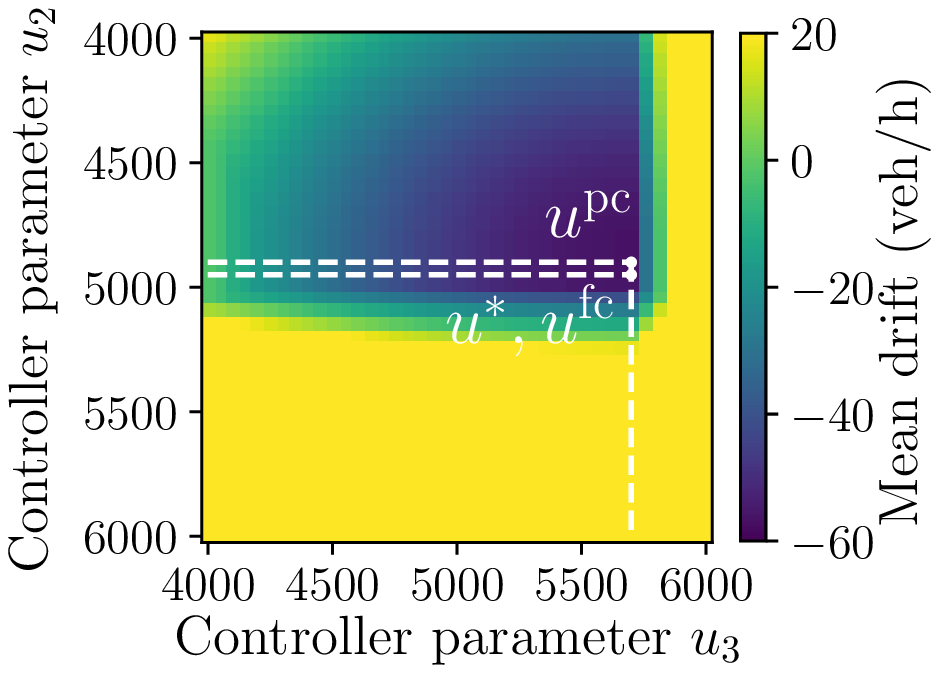}
\label{fig_Dbar2}
}
\subfigure[Simulated queue size.]{
\centering
\includegraphics[width=0.3\linewidth]{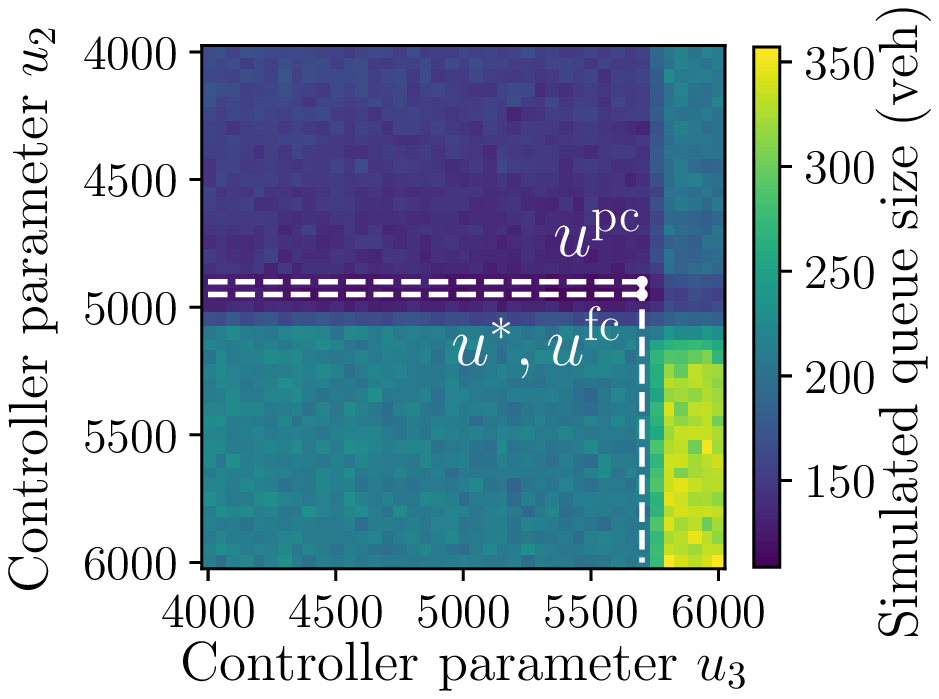}
\label{fig_qbar}
}
\caption{Comparison of (i) the solution $u^*$ to $\mathrm{P}_0$, (ii) the solution $u^{\mathrm{fc}}$ to $\mathrm{P}_3$, (iii) the solution $u^{\mathrm{pc}}$ to $\mathrm{P}_4^1$ and $\mathrm{P}_4^2$.}
\label{fig_Dbar_qbar}
\end{figure}

We first compute the mean drift \eqref{eq_psDs2} through the grid search of $u_2$ and $u_3$, as shown in Fig.~\ref{fig_Dbar2}. We found that $(u_2^{\mathrm{fc}}, u_3^{\mathrm{fc}})=(4950, 5700)$ minimizes the mean drift and hence is an optimal solution to $\mathrm{P}_3$. Then we applied the partially coordinated approach and obtained $(u_2^{\mathrm{pc}}, u_3^{\mathrm{pc}}) = (4900, 5700)$. Thus, in this example, the partially coordinated approach gave a \emph{near-optimal solution} to the coordinated control design problem $\mathrm{P}_3$. Finally, for different combinations of $(u_2, u_3)$, we used the numerical simulation to estimate the time-averaged buffer queue sizes $
\hat{Q}^{\mu}=\frac{1}{T}\sum_{t=1}^T (q_1^{\mu}(t)+q_2^{\mu}(t)+q_3^{\mu}(t))$, where $T$ denotes the number of simulation steps. We considered $T=10^6$ with a step size $\Delta_t=10$ seconds.
We found that $(u_2^*, u_3^*)=(4950,5700)$ minimizes the average simulated queue length and can be seen as the optimal solution to $\mathrm{P}_0$. Thus, the solution to $\mathrm{P}_3$ is approximately optimal to $\mathrm{P}_0$. The numerical simulation shows that $u^*$ and $u^{\mathrm{fc}}$ both result in $\hat{Q}^*=\hat{Q}^{\mathrm{fc}}=110$ veh and that $u^{\mathrm{pc}}$ leads to $\hat{Q}^{\mathrm{pc}}=120$ veh. This shows that the partially coordinated control design can yield the performance comparable to that of the fully coordinated control design. Our finding is important because the fully coordinated control design becomes notoriously difficult due to NP-hardness of such bi-level control problems. Our partially coordinated control offers a worthwhile approach to solving ramp metering for realistic highway in the face of capacity perturbations.

\begin{figure*}[h]
\centering
\subfigure[ALINEA.]{
\centering
\includegraphics[width=0.22\textwidth]{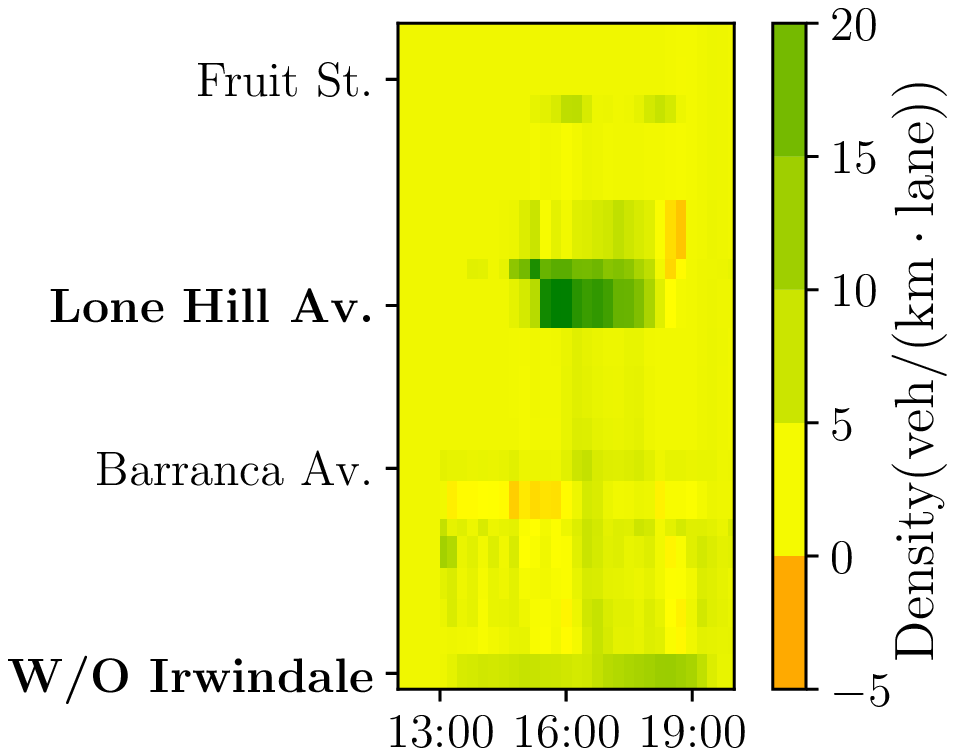}
\label{fig_case_alinea_speedmap}
}
\subfigure[METALINE.]{
\centering
\includegraphics[width=0.22\textwidth]{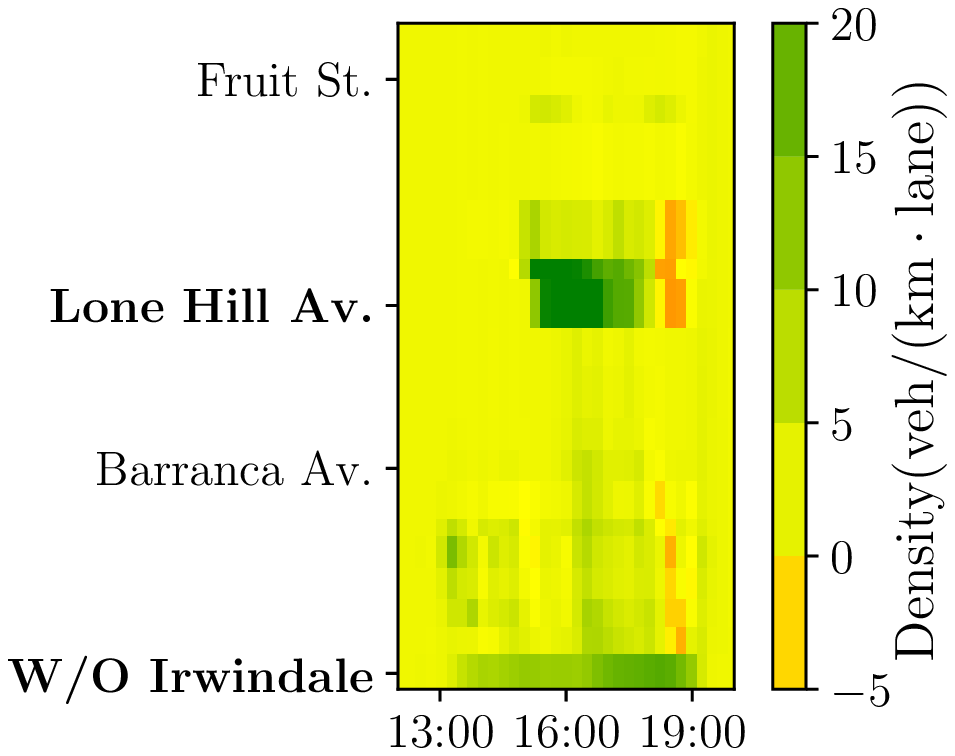}
\label{fig_case_loc_speedmap}
}
\subfigure[Localized.]{
\centering
\includegraphics[width=0.22\textwidth]{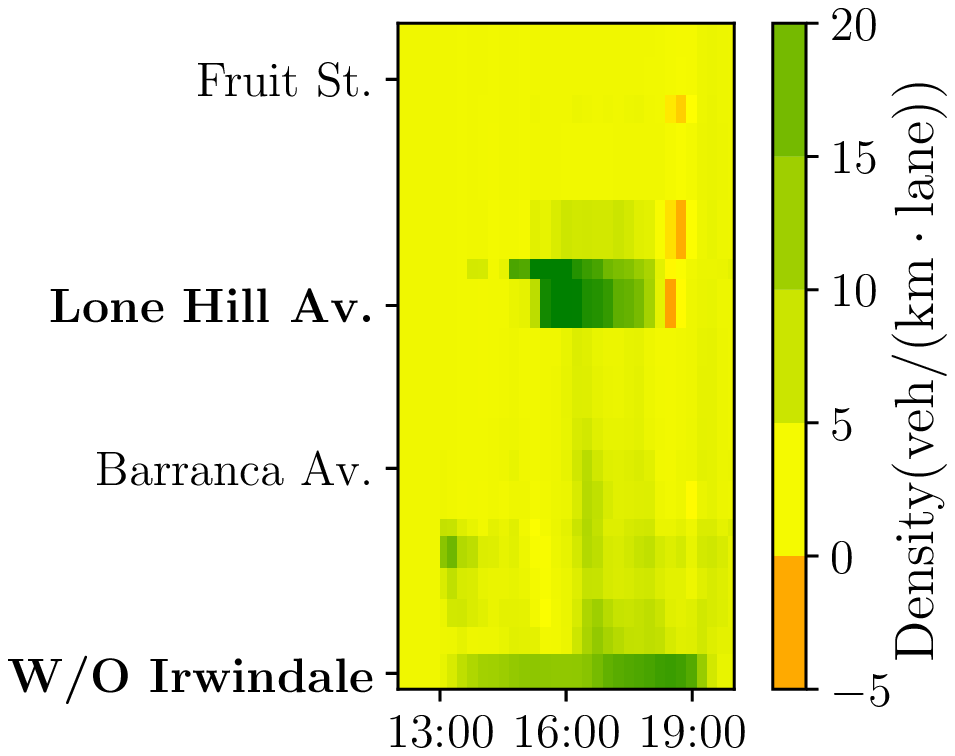}
\label{fig_case_metaline_denmap}
}
\subfigure[Partially coordinated.]{
\centering
\includegraphics[width=0.22\textwidth]{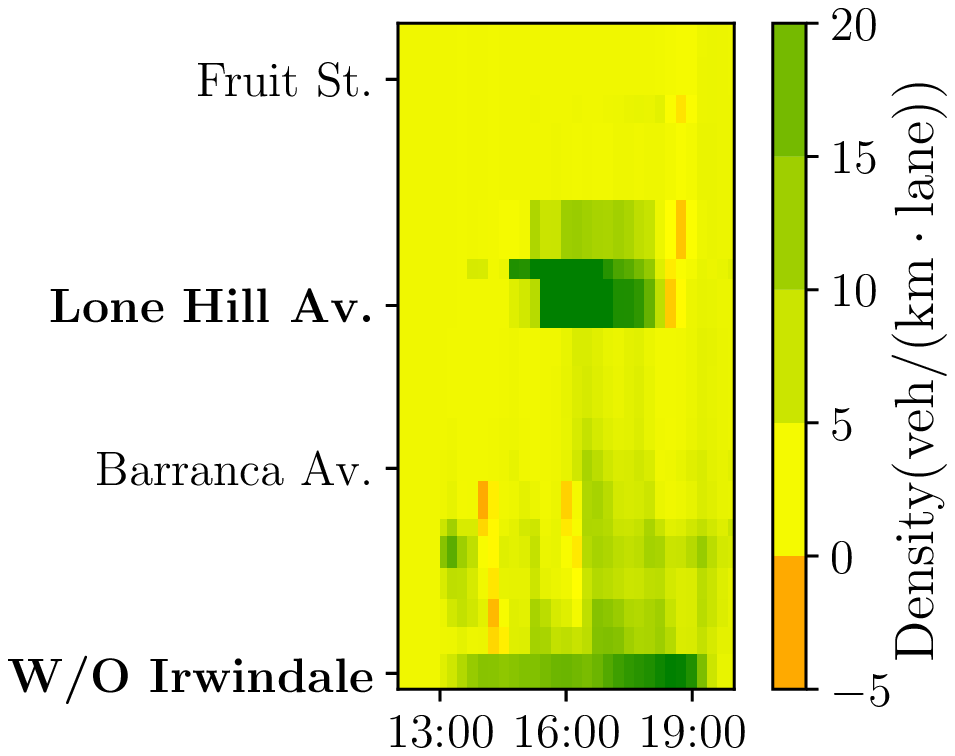}
\label{fig_case_cas_denmap}
}

\caption{Comparison of differences between traffic density maps with and without ramp metering.}
\label{fig_case_denmap}
\end{figure*}
\section{Case study: US Interstate 210}
\label{sec_simulate}

This section presents the design of ramp metering for the I-210E stretch that we introduced in Section I. Concretely, we applied the localized and partially coordinated control design based on the demand data and the calibrated SS-CTM. Although the proposed methods have been already illustrated with using numerical examples, this case study aims at evaluating and comparing them under more realistic conditions, in particular time-varying demand.

\begin{figure}[htbp]
    \centering
    \includegraphics[width=0.5\linewidth]{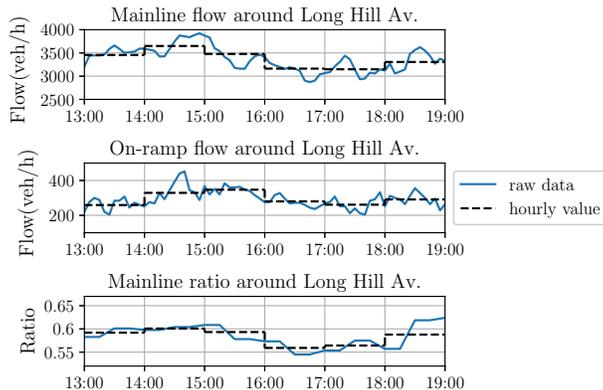}
    \caption{Hourly traffic demands around Long Hill Avenue.}
    \label{fig_hourly}
\end{figure}

We consider five scenarios. In the first scenario there is no ramp metering, and the remaining four adopt different control strategies. Recall from Fig.~\ref{fig_calibratedSpeedMap} that the traffic congestion occurred between 13:00 and 19:00; so we considered ramp metering turned on in that period. The second and the third scenarios correspond to the localized ALINEA and the coordinated METALINE, respectively. The METALINE controller is given by
\begin{equation}
\mu(t) = \mu(t-1) - K_P(n(t) -n(t-1)) - K_I(n(t)-n^c), \label{eq_metaline}
\end{equation}
where $K_P$ and $K_I$ are gain matrices, and $n^c\in\mathbb{R}^K$ are the nominal critical densities. We computed the control parameters $K_P$ and $K_I$ by solving a mixed-integer bilinear program that minimizes \emph{vehicle hours traveled} $H$ for the nominal CTM: 
\begin{equation}
    H = \Delta_t \Big(\sum_{k=1}^K \sum_{t=1}^T l_k n_k(t) + \sum_{k=1}^{K} \sum_{t=1}^T q_k(t) \Big), \label{eq_vht}
\end{equation}
where $\Delta_t$ denotes time step size. The optimal VHT between 13:00 and 19:00 attained 9829 veh$\cdot$hr under METALINE. In the last two scenarios, the localized and partially coordinated ramp metering were designed for each hour. We first computed hourly mainline flows, on-ramp flows and mainline ratios along the highway (see Fig.~\ref{fig_hourly}) and then used these hourly values for computing the control parameters.

The four ramp metering strategies were evaluated via numerical simulation. For SS-CTM the highway traffic was simulated from 12:00 to 21:00 with a time step size $\Delta_t=10$ seconds and each control strategy was tested for 1000 samples. 
While implementing ramp control, we also considered the capacity limits of on-ramp queues. In practice, ramp metering can result in long queues that spill from on-ramps and block street traffic. To avoid this, we restricted the maximal queue size at 40 veh/lane. Once the simulated queue size exceeded the threshold, the ramp metering would stop and on-ramp vehicles would enter the mainline at the maximum rate.

We consider two performance metrics, namely 1) time-averaged buffer queue length 
$$\hat{Q}=\frac{1}{T}\sum_{t=1}^T\sum_{k=1}^K q_k(t)$$
and 2) vehicle hours traveled \eqref{eq_vht}. The first is the control objective of the original control problem $\mathrm{P}_0$, and the second is commonly used for measuring the efficiency of highway operation, which includes mainline traveling time and queuing time of all vehicles.  Tables~\ref{tab_compareQ}-\ref{tab_compareVHT} summarize the evaluation of control strategies. 
Though the performance of our control design approach is marginal in the initial hours, they are effective in significantly reducing the queue length between 15:00-19:00. For example, compared with ALINEA (resp. METALINE), our localized control (resp. partially coordinated control) can shorten the queue length by 12.8\% (resp. 8.8\%) during 15:00-19:00. Overall, the localized and partially coordinated control designs reduce vehicle hours traveled by 8.3\% and 9.9\% respectively, while the classical ALINEA and METALINE decrease VHT by 5.1 \% and 6.2\%. This demonstrates that the highway system performance can be improved by considering stochastic capacity in designing ramp controllers. In this paper, we considered the affine controller \eqref{eq_controller} that relies on local measurements. Potentially higher gains can be achieved by considering more sophisticated control strategies based on more upstream or downstream measurements.

\begin{table}[htbp]
    \centering
    \scriptsize
    \caption{Time-averaged queue length (veh).}
    \begin{tabu}to\linewidth{X[4.3,c,m]X[1,c,m]X[1,c,m]X[1,c,m]X[1,c,m]X[1,c,m]X[1,c,m]X[1.1,c,m]}
    \toprule
      & 13:00-14:00 & 14:00-15:00 & 15:00-16:00 & 16:00-17:00 & 17:00-18:00 & 18:00-19:00 &  Mean \\
    \midrule
    No control & 545 & 802 & 771 & 856 & 1087 & 1148 & 868 \\
    ALINEA & 548 & 781 & 862 & 979 & 1041 & 909 & 853 \\
    METALINE & 555 & 741 & 811 & 909 & 983 & 875 & 812 \\
    Localized & 555 & 722 & 747 & 855 & 896 & 802 & 777 \\
    Partially coordinated & 555 & 722 & 764 & 855 & 900 & 742 & 756 \\
    \bottomrule
    \end{tabu}
    \label{tab_compareQ}
\end{table}

\begin{table}[htbp]
    \centering
    \scriptsize
    \caption{Vehicle hours traveled (veh$\cdot$hr).}
    \begin{tabu}to\linewidth{X[4.3,c,m]X[1,c,m]X[1,c,m]X[1,c,m]X[1,c,m]X[1,c,m]X[1,c,m]X[1.1,c,m]}
    \toprule
      & 13:00-14:00 & 14:00-15:00 & 15:00-16:00 & 16:00-17:00 & 17:00-18:00 & 18:00-19:00 & Sum \\
    \midrule
    No control & 1795 & 2177 & 2476 & 2862 & 3040 & 2578 & 14928 \\
    ALINEA & 1770 & 2112 & 2382 & 2719 & 2816 & 2366 & 14165 \\
    METALINE & 1764 & 2099 & 2348 & 2658 & 2787 & 2350 & 14006 \\
    Localized & 1761 & 2078 & 2320 & 2592 & 2691 & 2247 & 13689 \\
    Partially coordinated & 1760 & 2078 & 2282 & 2511 & 2619 & 2201 & 13451 \\
    \bottomrule
    \end{tabu}
    \label{tab_compareVHT}
\end{table}

Fig.~\ref{fig_case_denmap} visualizes the difference between traffic density maps with and without ramp metering. Each subfigure presents a density map under a particular control strategy minus the background density map without ramp control (see Fig.~2(c)). The positive values indicate that the traffic jam is alleviated by the ramp metering and the negative values imply that the congestion worsens. We note that ALINEA significantly mitigates the traffic jam around Lone Hill Avenue. Our localized ramp control further reduces the congestion upstream of Irwindale Avenue. This can be ascribed to the controller around Barranca Avenue that is aware of the capacity variation. Similarly, METALINE also reduces the upstream congestion by coordinated control, but our capacity-aware partial coordination further ameliorates congestion around Lone Hill Avenue and upstream of Barranca Avenue.
\section{Concluding remarks}
\label{sec_conclude}

In this paper, we developed a control-theoretic approach for determining resilient ramp meter configurations facing stochastic capacity perturbations. We modeled traffic dynamics under Markovian capacities and formulated the general coordinated control problem $\mathrm{P}_1$ based on the Foster-Lyapunov criterion. Then, we proposed two kinds of ramp metering design, namely localized and partially coordinated. These control designs were investigated numerically and analytically. The numerical study revealed that our control design approach can provide near-optimal solutions in terms of minimizing queue lengths and that the control efficacy is more significant as capacity fluctuations become more intense. We analytically showed that the localized and coordinated methods can produce the same decoupling impact, which prevents congestion from interfering upstream traffic. In general, the partially coordinated control considering interacting on-ramps is superior to the localized control, which was demonstrated by the case study of I-210E. This case study also showed our methods outperform the classical ALINEA and METALINE in the face of capacity perturbations. We note that the proposed designs of state-feedback ramp metering rely on calibrated Markov processes from historical data while ALINEA and MEATLINE do not. Thus, our design approach should not be viewed as a replacement to classical ramp metering. Rather, we suggest that ramp metering can be made more resilient to capacity perturbations using our approach. Future works involve investigating the approach to tuning PI controllers, such as ALINEA and METALINE, against capacity perturbations, and developing control designs adaptable to demand uncertainty.
\section*{Appendices}
\renewcommand{\thesubsection}{\Alph{subsection}}

\subsection{Calibration of SS-CTM}
\label{sub_calibration}

We consider a 18.1-km stretch of I-210E as a test case for the proposed ramp control; see Fig.~\ref{fig_map}. There are 11 on-ramps and 10 off-ramps. We divide the highway section into 17 cells based on the layout of on-/off-ramps and detectors. The calibration was built using Performance Measurement System (PeMS, see \cite{varaiya09}) data on Tuesday, March 26, 2019. This was one of the days when most detectors on mainline, on-ramps, and off-ramps were intact, and hence the PeMS data are reliable. 
\begin{table}[htbp]
    \centering
    \scriptsize
    \caption{Calibrated traffic flow parameters.}
    \begin{tabu}{ccccc}
    \toprule
    Cell $k$ & $v_k$ & $n_k^c$ & $w_k$ & $n_k^{\mathrm{jam}}$ \\
    & (km/h) & (veh$\cdot$lane/km) & (km/h) & (veh$\cdot$lane/km) \\
    \midrule
     1 & 99.0 & 85.7 & 13.7 & 598.2 \\
     2 & 99.1 & 72.3 & 11.7 & 664.0 \\
     3 & 99.7 & 71.9 & 18.4 & 460.2 \\
     4 & 99.3 & 74.5 & 12.7 & 664.0 \\
     5 & 99.5 & 69.9 & 11.7 & 664.0 \\
     6 & 96.2 & 73.4 & 11.1 & 664.0 \\
     7 & 97.4 & 74.1 & 11.5 & 664.0 \\
     8 & 97.1 & 74.4 & 20.5 & 400.5 \\
     9 & 96.8 & 77.4 & 24.8 & 363.0 \\
     10 & 96.2 & 77.9 & 19.7 & 424.0 \\
     11 & 99.4 & 64.8 & 21.1 & 664.0 \\
     12 & 102.7 &43.6  & 11.0 & 530.0 \\
     13 & 96.7 & 62.4  & 22.7& 339.5 \\
     14 & 95.7 & 57.8 & 11.5 & 498.0 \\
     15 & 97.7 & 53.7 & 12.5 & 498.0 \\
     16 & 102.3 &53.6 & 11.2 & 498.0 \\
     17 & 98.4 & 54.9 & 11.8 & 498.0 \\
    \bottomrule
    \end{tabu}
    \label{tab_calibratedTrafficParam}
\end{table}

The model calibration comprises two parts \cite{jin2017calibration}. The first part estimates nominal traffic flow parameters listed in Table~\ref{tab_calibratedTrafficParam}. Note that we restrict the jam density per lane not to exceed 166 veh/km. The second part involves the Markovian highway capacity. We first identified the capacity variation and found the cell 7 and 16 are subject to major perturbations; see Figs.~\ref{fig_cap_cell7}-\ref{fig_cap_cell16}. Since these two cells are far apart from each other, we assume their capacities fluctuate independently. For each cell, we considered that its capacity switched between two states, which is motivated by the two-capacity phenomenon \cite{banks1991two}. One can refer to \cite{jin2017calibration} for more details of calibrating Markov processes. The results are presented in Table~\ref{tab_markovprocess}, where $F_{k,1}$ (resp. $F_{k,2}$) is the nominal (resp. perturbed) capacity and $\lambda_{1,2}$ (resp. $\lambda_{2,1}$) is transition probability from $F_{k,1}$ (resp. $F_{k,2}$) to $F_{k,2}$ (resp. $F_{k,1}$) every hour. 
\begin{table}[htbp]
    \centering
    \scriptsize
    \caption{Calibrated Markovian capacity.}
    \begin{tabu}{ccc}
    \toprule
     & Cell 7 & Cell 16  \\
     & (Barranca Avenue) &  (Fruit Street) \\
    \midrule
    Period & 12:00-21:00 & 15:00-18:00 \\
    $F_{k,1}$ (veh/h) & 7224 & 5484 \\
    $F_{k,2}$ (veh/h) & 6670 & 4700 \\
    $\lambda_{1,2}$ (/h) & 0.6 & 0.6 \\
    $\lambda_{2,1}$ (/h) & 0.48 & 0.48 \\
    \bottomrule
    \end{tabu}
    \label{tab_markovprocess}
\end{table}


\subsection{Proof of Theorem~\ref{thm_distributed}}
\label{app_pf_thm1}
We construct a Lyapunov function and show that \eqref{eq_psDs} ensures the Foster-Lyapunov criterion for stability. First, note the invariant set $\mathcal{M}^{{\mathrm{lo}}}$ stated by Lemma~\ref{lmm_M}. The lemma is proved in Appendix~G.
\begin{lmm}\label{lmm_M}
Consider a set of the state $(q_1,q_2,n_1,n_2)$:
\begin{equation*}
    \mathcal M^{{\mathrm{lo}}}:=\bigcup_{i_1,i_2\in\{0,1\}}
    \mathcal{Q}_{1,i_1}^{{\mathrm{lo}}}\times\mathcal{Q}_{2,i_2}^{{\mathrm{lo}}}\times\mathcal{N}_{1,i_1}^{{\mathrm{lo}}}\times\mathcal{N}_{2,i_2}^{{\mathrm{lo}}},
\end{equation*}
where for $k=1,2$, $i_k$ indicates whether there is a queue in buffer $k$, and
\begin{align*}
    \mathcal{Q}_{k,i_k}^{{\mathrm{lo}}}:= \begin{cases}
    \{0\} & i_k=0 \\
    \mathbb{R}_{>0} & i_k=1
    \end{cases},~
    \mathcal{N}_{k,i_k}^{{\mathrm{lo}}}:= \begin{cases}
    [\underline{n}_k, \tilde{n}{_k^{{\mathrm{lo}}}}] & i_k=0 \\
    [{\uwave n}{_k^{{\mathrm{lo}}}}, \tilde{n}{_k^{{\mathrm{lo}}}}] & i_k=1
    \end{cases}.
\end{align*}
The set $\mathcal M^{{\mathrm{lo}}}$ is an invariant set for the two-cell highway section.
\end{lmm}

Then we consider the Lyapunov function $V:\mathcal{S}\times\mathcal{M}^{\mathrm{lo}}\to\mathbb{R}_{\geq0}$ as follows:
\begin{align}
    &V^{\mathrm{lo}}(\phi) :=  q_1(\frac{1}{2} q_1 + \frac{1}{2}\beta_1q_2 + l_1n_1 + \beta_1l_2\int_{\underline n_2}^{n_2} \rho_2(\xi) \mathrm{d}\xi+b_{1,s}) \nonumber \\
    & + q_{2}(\frac{1}{2}\beta_1q_1+\frac{1}{2}q_2 + \beta_1l_1n_1 +  l_2\int_{\underline n_2}^{n_2} \rho_2(\xi) \mathrm{d}\xi + b_{2,s}), \nonumber
\end{align}
where the non-negative parameters $b_{k,s}$, $s\in\mathcal S, k\in\{1, 2\}$, are solutions to
$D^{{\mathrm{lo}}}_{k,s}(\mathcal{E}_k^{{\mathrm{lo}}}) + \sum_{s'\in\mathcal{S}}\lambda_{s,s'}(b_{k,s'}-b_{k,s}) = \sum_{s'\in\mathcal{S}} p_s D^{{\mathrm{lo}}}_{k,s'}(\mathcal{E}_k^{{\mathrm{lo}}})$. The existence of $b_{k,s}$ is guaranteed by Lemma~\ref{lmm_bs} proved in Appendix~H.
\begin{lmm}\label{lmm_bs}
The following system of equations, given a set of $z_{k,s}$ with $s\in\mathcal{S}$ and $k=\{1,\cdots,K\}$,
\begin{equation*}
    z_{k,s} + \sum_{s'\in\mathcal{S}}\lambda_{s,s'}(b_{k,s'}-b_{k,s}) = \sum_{s'\in\mathcal{S}}p_sz_{k,s'}
\end{equation*}
has a non-negative solution for $b_{k,s}$.
\end{lmm}

Applying the infinitesimal generator, we have
\begin{align*}
    \mathscr LV^{\mathrm{lo}}(\phi) \leq& \Big(D_1^{\mathrm{lo}}(\phi) + \sum_{s'\in\mathcal S}\lambda_{s,s'}(b_{1,s'}-b_{1,s})\Big)q_1 \\
    & + \Big(D_2^{\mathrm{lo}}(\phi) + \sum_{s'\in\mathcal S}\lambda_{s,s'}(b_{2,s'}-b_{2,s})\Big) q_2  + d_0,
\end{align*}
where the existence of $d_0<\infty$ is guaranteed by the finite $G_k^{\mathrm{lo}}$, $H_k^{\mathrm{lo}}$ and $n_k$ for $k=1,2$.

The rest of this proof is devoted to verifying the drift condition
\begin{equation}
    \mathscr L V^{\mathrm{lo}}(\phi)\le-c^{{\mathrm{lo}}}(q_1+q_2)+d_0,~\forall \phi\in\mathcal{S}\times\mathcal{M}^{{\mathrm{lo}}} \label{eq_thm1_pf_1}
\end{equation}
for some $c^{{\mathrm{lo}}}>0$. We consider $\mathcal M^{{\mathrm{lo}}}$ by discussing whether there are queues in both buffer 1 and 2. Let $\mathcal{M}_{i_1,i_2}^{{\mathrm{lo}}}:=\mathcal{Q}_{1,i_1}^{{\mathrm{lo}}}\times\mathcal{Q}_{2,i_2}^{{\mathrm{lo}}}\times\mathcal{N}_{1,i_1}^{{\mathrm{lo}}}\times\mathcal{N}_{2,i_2}^{{\mathrm{lo}}}$. 

Over the $\mathcal{S}\times\mathcal{M}_{0,0}^{{\mathrm{lo}}}$, we have $q_1=q_2=0$. Clearly, \eqref{eq_thm1_pf_1} holds. Over the set  $\mathcal{S}\times(\mathcal{M}_{1,0}^{{\mathrm{lo}}}\cup\mathcal{M}_{0,1}^{{\mathrm{lo}}}\cup\mathcal{M}_{1,1}^{{\mathrm{lo}}})$, we have
$D^{{\mathrm{lo}}}_{1,s}(\mathcal{M}_{1,0}^{{\mathrm{lo}}}\cup\mathcal{M}_{0,1}^{{\mathrm{lo}}}\cup\mathcal{M}_{1,1}^{{\mathrm{lo}}})q_1 = D^{{\mathrm{lo}}}_{1,s}(\mathcal M^{{\mathrm{lo}}}_{1,0}\cup\mathcal M^{{\mathrm{lo}}}_{1,1})q_1$
by noting $q_1=0$ for $(q,n)\in\mathcal{M}_{0,1}^{{\mathrm{lo}}}$. Next, we show $D^{{\mathrm{lo}}}_{1,s}(\mathcal M^{{\mathrm{lo}}}_{0,1}\cup\mathcal M^{{\mathrm{lo}}}_{1,1})=D^{{\mathrm{lo}}}_{1,s}(\mathcal E^{{\mathrm{lo}}}_{1})$. Substituting \eqref{eq_G1}-\eqref{eq_H} into $D_1^{\mathrm{lo}}(\phi)$ yields
\begin{align*}
    D^{{\mathrm{lo}}}_{1}(\phi) =& \alpha_1 + \beta_1\alpha_2 - (1-\beta_1^2 \rho_2(n)) f_1^{{\mathrm{lo}}}(\phi) \\
    & -  \beta_1(1-\rho_2(n)) r_2^{{\mathrm{lo}}}(q_2,n_2) - \beta_1\rho_2(n) f_2^{{\mathrm{lo}}}(\phi).
\end{align*}
Note that $q_1$ is cancelled and that $D^{{\mathrm{lo}}}_1(\phi)$ is monotonically decreasing in $n_1$. So we conclude
$$D^{{\mathrm{lo}}}_{1}(\phi)q_1\leq D^{{\mathrm{lo}}}_{1,s}(\mathcal E^{{\mathrm{lo}}}_{1})q_1,~\forall \phi\in\mathcal{S}\times(\mathcal{M}_{1,0}^{{\mathrm{lo}}}\cup\mathcal{M}_{0,1}^{{\mathrm{lo}}}\cup\mathcal{M}_{1,1}^{{\mathrm{lo}}}).$$
Similarly, we can show 
$$D^{{\mathrm{lo}}}_{2}(\phi)q_2\leq D^{{\mathrm{lo}}}_{2,s}(\mathcal E^{{\mathrm{lo}}}_{2})q_2,~\forall \phi\in\mathcal{S}\times(\mathcal{M}_{1,0}^{{\mathrm{lo}}}\cup\mathcal{M}_{0,1}^{{\mathrm{lo}}}\cup\mathcal{M}_{1,1}^{{\mathrm{lo}}}).$$
It follows $\mathscr LV^{\mathrm{lo}}(\phi)
    \leq \bar{D}^{\mathrm{lo}}(q_1+q_2) + d_0$ for any $\phi\in\mathcal{S}\times(\mathcal{M}_{1,0}^{{\mathrm{lo}}}\cup\mathcal{M}_{0,1}^{{\mathrm{lo}}}\cup\mathcal{M}_{1,1}^{{\mathrm{lo}}})$. Hence, there exists $c^{{\mathrm{lo}}}:=-\bar{D}^{{\mathrm{lo}}} > 0$ given \eqref{eq_drift}.

Combining the results above, we conclude that the drift condition \eqref{eq_drift} is achieved with $c^{{\mathrm{lo}}}=-\bar{D}^{{\mathrm{lo}}}$ and $d=d_0$, which implies the stability. \qed

\subsection{Proof of Theorem~\ref{thm_centralized}}

First, we have the invariant set $\mathcal{M}^{{\mathrm{fc}}}$ stated by the following lemma:
\begin{lmm} \label{lmm_M_co}
The set
\begin{equation}
    \mathcal{M}^{{\mathrm{fc}}}:=\bigcup_{
    i_j\in\{0,1\},1\leq j\leq K} \prod_{j=1}^K
    \mathcal{Q}_{j,i_j}^{{\mathrm{fc}}}\times
    \prod_{j=1}^K  \mathcal{N}_{j,i_j}^{{\mathrm{fc}}}
    \label{eq_M_fc}
\end{equation}
is an invariant set for the $K$-cell highway section.
\end{lmm}

The lemma is proved in Appendix~I. Note that the definition yields $\mathcal{E}_k^{\mathrm{fc}}=\{(q,n)\in\mathcal{M}^{\mathrm{fc}}|q_k=1\}$.

Then we consider the Lyapunov function $V^{\mathrm{fc}}:\mathcal{S}\times\mathcal{M}^{\mathrm{fc}}\to\mathbb{R}_{\geq0}$.
\begin{align*}
    V^{\mathrm{fc}}(\phi) =& \sum_{k=1}^K
    q_k\bigg(\sum_{j=1}^{k-1}\gamma_{j,k}\Big(\frac{1}{2}q_j + l_j\int_{\underline{n}_j}^{n_j} \rho_j(\xi) \mathrm{d}\xi  \Big) \\
    & + \sum_{j=k}^K\gamma_{k,j} \Big(\frac{1}{2}q_j+l_j\int_{\underline{n}_j}^{n_j} \rho_j(\xi) \mathrm{d}\xi\Big) +  b_{k,s} \bigg).
\end{align*}
The rest of the proof is similar to that of Theorem~\ref{thm_distributed}. \qed

\subsection{Proof of Theorem~\ref{thm_cascaded}}
First consider an invariant set $\mathcal{M}^{\mathrm{pc}}$ analogous to \eqref{eq_M_fc} but under the partially coordinated control policy $\mu^{\mathrm{pc}}$. Then the proof uses the Lyapunov function $V^{\mathrm{pc}}:\mathcal{S}\times\mathcal{M}^{\mathrm{pc}}\to\mathbb{R}_{\geq0}$:
\begin{align*}
    V^{\mathrm{pc}}(\phi) =& \sum_{k=1}^K 
    q_k\bigg(\sum_{j=1}^{k-1}\gamma_{j,k}\Big(\frac{1}{2}q_j + l_jn_j \Big) \\
    & + \sum_{j=k}^K\gamma_{k,j} \Big(\frac{1}{2}q_j+l_j\int_{\underline{n}_j}^{n_j} \rho_j(\xi) \mathrm{d}\xi\Big) +  b_{k,s} \bigg).
\end{align*}
The remaining proof is similar to that of Theorem~\ref{thm_distributed}.

\subsection{Proof of Corollary~1}

We first consider the necessity. If the two-cell highway section is stable, we must conclude the time-averaged inflow is less than the time-averaged outflow, which indicates \eqref{eq_pr1_conc1}-\eqref{eq_pr1_conc2}. 

Then, we show the sufficiency. Consider the invariant set $\tilde{\mathcal M}^{{\mathrm{lo}}}$
stated by Lemma~\ref{lmm_tighterM}.
\begin{lmm}\label{lmm_tighterM}
The set 
$$\tilde{\mathcal M}^{{\mathrm{lo}}}:=\bigcup_{i_1,i_2\in\{0,1\}}
\tilde{\mathcal{Q}}_{1,i_1,i_2}^{{\mathrm{lo}}}\times\tilde{\mathcal{Q}}_{2,i_1,i_2}^{{\mathrm{lo}}}\times\tilde{\mathcal{N}}_{1,i_1,i_2}^{{\mathrm{lo}}}\times\tilde{\mathcal{N}}_{2,i_1,i_2}^{{\mathrm{lo}}} $$ 
is an invariant set for the two-cell highway section satisfying \eqref{eq_decoupling} and \eqref{eq_pr1_con1}-\eqref{eq_pr1_con3},
where for $k=1,2$,
\begin{equation*}
    \tilde{\mathcal{Q}}_{k,i_1,i_2}^{{\mathrm{lo}}}:= \begin{cases}
    \{0\} & i_k=0, \\
    \mathbb{R}_{>0} & i_k=1,
    \end{cases}, \tilde{\mathcal{N}}_{k,i_1,i_2}^{{\mathrm{lo}}}:= \begin{cases}
    [\underline{n}_k, \bar{n}_k] & i_j=0, \forall 1\leq j\leq k, \\
    [n_k^c, \bar{n}_k] & \mathrm{otherwise}.
    \end{cases}
\end{equation*}
\end{lmm}

The lemma is proved in Appendix~J. Note that $\tilde{\mathcal M}^{{\mathrm{lo}}}$ is in general tighter than $\mathcal M^{{\mathrm{lo}}}$ and leads to tighter $\tilde{\mathcal{E}}_1^{{\mathrm{lo}}} :=\{0\}\times\{n_1^c\}\times[n_2^c, \bar{n}_2]$ and $\tilde{\mathcal{E}}_2^{{\mathrm{lo}}} :=\{1\}\times\{\underline n_1\} \times [n_2^c, \bar{n}_2]$.

According to Theorem~\ref{thm_distributed}, we only need to consider $D_{1,s}^{{\mathrm{lo}}}(\tilde{\mathcal{E}}_1^{{\mathrm{lo}}})$ and $D_{2,s}^{{\mathrm{lo}}}(\tilde{\mathcal{E}}_2^{{\mathrm{lo}}})$. Considering $\rho_2(n)=1$ in case of the decoupled cells 1 and 2, we have
\begin{align*}
     D_{1,s}^{{\mathrm{lo}}}(\tilde{\mathcal{E}}_1^{{\mathrm{lo}}}) =& \max_{(q,n)\in\tilde{\mathcal{E}}_1^{{\mathrm{lo}}}}(1-\beta_1^2)(\alpha_1-\min\{v_1n_1,F_{1,s}\}) \\
     & + (\beta_1\alpha_1+\alpha_2-\min\{v_2n_2,F_{2,s}\}) \\
     =& (1-\beta_1^2)(\alpha_1-F_{1,s}) + (\beta_1\alpha_1+\alpha_2-F_{2,s}).
\end{align*}
Similarly, we obtain $D_{2,s}^{{\mathrm{lo}}}(\tilde{\mathcal{E}}_2^{{\mathrm{lo}}}) = \beta_1\alpha_1 + \alpha_2 - F_{2,s}$.
So given \eqref{eq_pr1_conc1}-\eqref{eq_pr1_conc2}, the stability condition \eqref{eq_psDs} is satisfied. \qed

\subsection{Proof of Corollary~\ref{prp_equi}}
We first consider the invariant set $\tilde{\mathcal M}^{\mathrm{pc}}$ stated in Lemma~\ref{lmm_tighterM_K} proved in Appendix~K. 
\begin{lmm}\label{lmm_tighterM_K}
The set 
\begin{equation*}
    \tilde{\mathcal{M}}^{\mathrm{pc}}:=\bigcup_{
    i_j\in\{0,1\},1\leq j\leq K} \prod_{j=1}^K
    \tilde{\mathcal{Q}}_{j,i_1,\cdots,i_K}^{\mathrm{pc}}\times
    \prod_{j=1}^K  \tilde{\mathcal{N}}_{j,i_1,\cdots,i_K}^{\mathrm{pc}}
\end{equation*}
is an invariant set for the $K$-cell highway section satisfying \eqref{eq_thm3_con} and \eqref{eq_pr2_con1}-\eqref{eq_pr2_con3}, where for $k=1,\cdots,K$,
\begin{subequations}
    \begin{align}
    &\tilde{\mathcal{Q}}_{k,i_1,\cdots,i_K}^{\mathrm{pc}}:= \begin{cases}
    \{0\} & i_k=0, \\
    \mathbb{R}_{>0} & i_k=1,
    \end{cases} \label{eq_invariantset_ca_1}\\
    &\tilde{\mathcal{N}}_{k,i_1,\cdots,i_K}^{\mathrm{pc}}:= \begin{cases}
    [\underline{n}_k, \bar{n}_k] & i_j=0, ~\forall 1\leq j \leq k, \\
    [n_k^c, \bar{n}_k] & \mathrm{otherwise}.
    \end{cases} \label{eq_invariantset_ca_2}
\end{align}
\end{subequations}
\end{lmm}

First, we show that $\mu^{\mathrm{pc}}$ is an optimal solution to the program $\mathrm{P}_3$ with $\rho_j\equiv1$, $j=1,2,\cdots,K$. For any $\mu^{\mathrm{fc}}$, we must have
\begin{align*}
    &\bar{D}^{{\mathrm{fc}}} \geq  \sum_{s\in\mathcal{S}}p_s\bigg((1-\beta_k^2)\Big(\sum_{j=1}^k \gamma_{1,j}\alpha_j - F_{k,s}\Big) \\
    & + \sum_{j=k+1}^K\gamma_{k,j}(1-\beta_j^2)\Big(\sum_{\ell=1}^j \gamma_{\ell,j}\alpha_j - F_{j,s}\Big)\bigg), ~\forall k=1,\cdots,K.
\end{align*}
Note that the control policy $\mu^{\mathrm{pc}}$ with the invariant set $\tilde{\mathcal{M}}^{\mathrm{pc}}$ achieves the equality. Thus we can conclude that $\mu^{\mathrm{pc}}$ is the optimal solution to $\mathrm{P}_3$.

Similarly, we can prove that $\mu^{\mathrm{pc}}$ can be obtained from the localized control design. \qed

\subsection{Proof of Lemma~\ref{lmm_M}}

We prove this result by showing that for each mode $s\in\mathcal{S}$, the vector field $(G,H)$ points towards to the interior on the boundary of $\mathcal M^{{\mathrm{lo}}}$. The boundary is the union of the following six sets $\mathcal{B}_i^{{\mathrm{lo}}}$, $i=1,\cdots,6$. We consider them one by one:

\begin{itemize}
    \item $\mathcal{B}_1^{{\mathrm{lo}}}:=\{(q,n)\in\mathcal M^{{\mathrm{lo}}}:q_1=0,n_1=\underline n_1\}$. Over this set, we have
    \begin{align*}
    H^{{\mathrm{lo}}}_1(\phi)
    & \geq (\min\{\alpha_1, w_1(n_1^{\mathrm{jam}} - \underline{n}_1)\} - v_1\underline{n}_1 ) / l_1 \\ 
    & \overset{\eqref{eq_fun}}{\geq} (\min\{\alpha_1, F_1^{\max}\}-v_1\underline n_1 ) / l_1 = 0.
    \end{align*}

    \item $\mathcal{B}_2^{{\mathrm{lo}}}:=\{(q,n)\in\mathcal M^{{\mathrm{lo}}}:q_2=0,n_2=\underline n_2\}$.
    Over this set, we have
    \begin{align*}
    H^{{\mathrm{lo}}}_2(\phi)
    \ge& (\min\{\beta_1 v_1 \underline n_1 + r_2^{{\mathrm{lo}}}(q_2, n_2), \\
    &\beta_1 F^\min_1 + r_2^{{\mathrm{lo}}}(q_2,n_2), F_2^{\max}\} - v_2\underline n_2 ) / l_2 \\
    =& 0.
    \end{align*}
    Note that we assume $r_2^{{\mathrm{lo}}}(q_2,n_2) \ge \alpha_2$ given $q_2=0$ and $n_2 = \underline n_2$; otherwise the controller is unreasonable and must lead to the instability. 

    \item $\mathcal{B}_3^{{\mathrm{lo}}}:=\{(q,n)\in\mathcal{M}^{{\mathrm{lo}}}:q_1>0,n_1=\uwave n{_1}\}$. The proof is similar to that of $\mathcal{B}_1^{{\mathrm{lo}}}$.
    
    \item $\mathcal{B}_4^{{\mathrm{lo}}}:=\{(q,n)\in\mathcal M^{{\mathrm{lo}}}:q_2>0,n_2=\uwave n{_2^{{\mathrm{lo}}}}\}$. The proof is similar to that of $\mathcal{B}_2^{{\mathrm{lo}}}$.
    
    \item $\mathcal{B}_5^{{\mathrm{lo}}}:=\{(q,n)\in\mathcal M^{{\mathrm{lo}}}:n_1=\tilde n_1^{{\mathrm{lo}}}\}$. Over this set, we have
    \begin{align*}
        H_1^{{\mathrm{lo}}}(\phi) 
        \le& (w_1(n_1^{\mathrm{jam}} - \tilde{n}_1^{{\mathrm{lo}}}) - \min\{F_1^{\min}, \\
        &(w_2(n_2^{\mathrm{jam}} - n_2) - r_2^{{\mathrm{lo}}}(1, n_2))/\beta_1 \} )/l_1 = 0.
    \end{align*}
    
    \item $\mathcal{B}_6^{{\mathrm{lo}}}:=\{(q,n)\in\mathcal M^{{\mathrm{lo}}}:n_2=\bar n_2\}$.
    The proof is similar to that of $\mathcal{B}_5^{{\mathrm{lo}}}$. \qed
\end{itemize}

\subsection{Proof of Lemma~\ref{lmm_bs}}
\label{app_bs}
First, we define the matrices
\begin{equation*}
    \Lambda:=\left[\begin{array}{ccc}
        -\sum_{i\neq1}\lambda_{1,i}  & \cdots & \lambda_{1,m} \\
        \vdots & \ddots & \vdots\\
        \lambda_{m,1} & \cdots & -\sum_{i\neq m}\lambda_{m,i}
    \end{array}\right], P:=\left[\begin{array}{ccc}
    p_1 &    & \\
        &    \ddots & \\
    p_1 & \cdots & p_m
    \end{array}\right],
\end{equation*}
where $\lambda_{m,i}$ denotes the transition rate from mode $m$ to mode $i$ and $p_i$ denotes the steady-state probability of mode $i$. Besides, let $p=[p_1,\cdots,p_m]^{\mathrm{T}}$, $b_k=[b_{k,1},\cdots,b_{k,m}]^{\mathrm{T}}$ and $z_k=[z_{k,1},\cdots,z_{k,m}]^{\mathrm{T}}$. To show the existence of a solution, note that the system of equations is equivalent to
$\Lambda b_k= [p^{\mathrm{T}}z_k-z_{k,1}, \cdots, p^{\mathrm{T}}z_k-z_{k,m}]^{\mathrm{T}}$. Since the discrete state process is ergodic, the rank of the matrix $\Lambda$ is $m-1$.
Scaling each row $i$ with $p_i$ and adding the scaled rows $1,2,\cdots,m-1$ to row $m$, we obtain $P\Lambda b_k = [p_1(p^{\mathrm{T}} z_k-z_{k,1}),\cdots,(\sum_{i=1}^mp_i)p^{\mathrm{T}}z_k - \sum_{i=1}^mp_iz_{k,i}]^{\mathrm{T}}$. Noting $\sum_{i\neq j}p_i\lambda_{i,j}-p_j\sum_{i\neq j}\lambda_{j,i}=0$ for $j=1,\cdots,m$ and $
(\sum_{i=1}^mp_i)p^{\mathrm{T}}z_k-\sum_{i=1}^mp_iz_{k,i}=p^{\mathrm{T}}z_k-p^{\mathrm{T}}z_k=0$, we conclude that the rank of the augmented coefficient matrix of the system of linear equations above is also $m-1$, equal to the rank of the coefficient matrix. Therefore, the system of equations above must have a solution $\hat{b}_k$. 
If $\hat{b}_k$ has some negative element $\hat{b}_{k,i} < 0$, we can obtain non-negative solution $\tilde{b}_k:=\hat{b}_k+|\min_i \hat{b}_{k,i}|$. \qed

\subsection{Proof of Lemma~\ref{lmm_M_co}}
\label{app_lmm4}
We show that for each mode $s\in\mathcal{S}$, the vector field $(G,H)$ points towards to the interior on the boundary of $\mathcal M^{{\mathrm{fc}}}$. The boundary of $\mathcal{M}^{{\mathrm{fc}}}$, denoted by $\mathcal{B}^{{\mathrm{fc}}}$, is the union of $3K$ sets:
$\mathcal{B}^{{\mathrm{fc}}}:=\bigcup_{i\in\{1,2,3\},k\in\{1,\cdots,K\}} \mathcal{B}^{{\mathrm{fc}}}_{i,k}$, where $\mathcal{B}^{{\mathrm{fc}}}_{1,k}:=\{(q,n)\in\mathcal M^{{\mathrm{fc}}}:q_k=0,n_k=\underline n_k\}$, $\mathcal{B}^{{\mathrm{fc}}}_{2,k}:=\{(q,n)\in\mathcal M^{{\mathrm{fc}}}:q_k>0,n_k=\uwave n{_k^{{\mathrm{fc}}}}\}$ and $\mathcal{B}^{{\mathrm{fc}}}_{3,k}:=\{(q,n)\in\mathcal M^{{\mathrm{fc}}}:n_k=\bar n{_k^{{\mathrm{fc}}}}\}$. The proofs for $\mathcal{B}^{{\mathrm{fc}}}_{i,k}$ are analogous to those for Lemma~\ref{lmm_M}. \qed

\subsection{Proof of Lemma~\ref{lmm_tighterM}}
\label{app_lmm3}
We first consider $\tilde{\mathcal{M}}_{1,0}^{{\mathrm{lo}}}$. Due to \eqref{eq_pr1_con1}, we have $\uwave{n}{_1}=\min\{U_1, F_1^{\max}\}/v_1 = n_1^c$. For $q_1 > 0$, $n_1 \ge n_1^c$, $q_2 = 0$ and $n_2 = n_2^c$, we have
\begin{align*}
    H_2^{{\mathrm{lo}}}(\phi)
    \ge& (\min\{\beta_1 F_1^{\min}+r_2^{{\mathrm{lo}}}(0, n_2^c), F_2^{\max}\} - F_2^{\max} ) / l_2 \geq 0.
\end{align*}
Hence, $n_2^c$ is one lower boundary. The proof is analogous for the boundary of $\tilde{\mathcal M}^{{\mathrm{lo}}}_{0,1}$ and $\tilde{\mathcal M}^{{\mathrm{lo}}}_{1,1}$. \qed

\subsection{Proof of Lemma~\ref{lmm_tighterM_K}}
\label{app_lmm5}
We show \eqref{eq_invariantset_ca_2} by induction. First, we have $\tilde{\mathcal{N}}_{1,0,i_2,\cdots,i_K}^{\mathrm{pc}}=[\underline{n}_1,\bar{n}_1]$ and $\tilde{\mathcal{N}}_{1,1,i_2,\cdots,i_K}^{\mathrm{pc}}=[n_1^c,\bar{n}_1]$ due to \eqref{eq_pr2_con1} for any $i_2,\cdots,i_K\in\{0,1\}$.  Then we assume \eqref{eq_invariantset_ca_2} for $k=i$. The proof for the case $k=i+1$ is similar to that for Lemma~\ref{lmm_tighterM}. \qed

\bibliographystyle{unsrt}

\bibliography{Bibliography}

\end{document}